\documentclass{article}[12pt]

\usepackage{xcolor,graphicx,float}
\usepackage{amsfonts,latexsym,amssymb}
\usepackage{mathrsfs,amsmath}
\DeclareFontFamily{U}{mathx}{\hyphenchar\font45}
\DeclareFontShape{U}{mathx}{m}{n}{
      <5> <6> <7> <8> <9> <10>
      <10.95> <12> <14.4> <17.28> <20.74> <24.88>
      mathx10
      }{}
\DeclareSymbolFont{mathx}{U}{mathx}{m}{n}
\DeclareFontSubstitution{U}{mathx}{m}{n}
\DeclareMathSymbol{\bigtimes}{1}{mathx}{"91}

\usepackage{amsfonts}

\newtheorem{Theorem}{Theorem}

\newtheorem{Proposition}[Theorem]{Proposition}

\newtheorem{Remark}[Theorem]{Remark}

\newcommand{\cvd}{\hfill{$\sqcap\!\!\!\!\sqcup$}}
\begin{document}

\title{A hybrid differential game with switching thermostatic-type dynamics and cost}

\author{Fabio Bagagiolo\footnote{Department of Mathematics, University of Trento, Italy  email {fabio.bagagiolo@unitn.it}}, Rosario Maggistro\footnote{Department of Management, Ca' Foscari University of Venice, Italy email {rosario.maggistro@unive.it}} and Marta Zoppello\footnote{Department of Mathematical Sciences, Politecnico di Torino, Italy email {marta.zoppello@polito.it}}}

\date{}

\maketitle

\begin{abstract}
In this paper, we consider an infinite horizon zero-sum differential game where the dynamics of each player and the running cost depend on the evolution of some discrete (switching) variables. In particular, such switching variables evolve according to the switching law of a so-called thermostatic delayed relay, applied to the players' states. We first address the problem of the continuity of both lower and upper value function. Then, by a suitable representation of the problem as a coupling of several exit-time differential games, we characterize those value functions as, respectively, the unique solution of a coupling of several Dirichlet problems for Hamilton-Jacobi-Isaacs equations. The concept of viscosity solutions and a suitable definition of boundary conditions in the viscosity sense are used in the paper.


\end{abstract}
\section{Introduction}

We consider an infinite horizon zero-sum differential game where both players, in their decoupled dynamics, as well as the running cost, are also affected by some switching variable whose switching evolution is described by a so-called thermostatic-type switching rule, subject to the evolution of the players' continuous state-variables. 

More precisely, we consider two decoupled dynamics for the two players with, respectively, state-variables $X\in\mathbb{R}^n$ and $Y\in\mathbb{R}^m$, as
\begin{equation}
\small
\label{eq:thermostat_systems_introduction}
\left\{
\begin{array}{ll}
X'(t)=f(X(t),W(t),\alpha(t)),\ t>0\\
W(t)=h[X, w](t)\\
X(0)=x, \ W(0)=w
\end{array}
\right.,
\ \ \ 
\left\{
\begin{array}{ll}
Y'(t)=g(Y(t),Z(t),\beta(t)),\ t>0\\
Z(t)=h[Y,z](t)\\
Y(0)=y, \ Z(0)=z
\end{array}
\right.
\end{equation}

\noindent
where $\alpha\in{\cal A},\beta\in{\cal B}$ are the measurable controls, and $W,Z\in\{-1,1\}$ are the switching variables, whose state-dependent switching rules are represented by the second lines of the systems \eqref{eq:thermostat_systems_introduction}. In particular, the switch from $1$ to $-1$ and the switch from $-1$ to $1$ occur at two different thresholds that must be reached by the continuous state variable  (see Figure \ref{fig:relay}, and in general Section \ref{sec:assumptions} for more precise details). Moreover, player $X$ wants to minimize, whereas player $Y$ wants to maximize, a discounted infinite horizon cost of the form

\[
J(x,y,w,z,\alpha,\beta)=\int_0^{+\infty}e^{-\lambda t}\ell(X(t),Y(t),W(t),Z(t),\alpha(t),\beta(t))dt
\]

\noindent
where $w,z\in\{-1,1\}$ are the initial states of the switching variables.

As usual, following Elliot-Kalton \cite{ellkal}, we define the lower and the upper value functions respectively as

\begin{equation}
\label{eq:values_introduction}
\begin{array}{ll}
\displaystyle
\underline V(x,y,w,z)=\inf_{\gamma\in\Gamma}\sup_{\beta\in{\cal B}} J(x,y,w,z,\gamma[\beta],\beta),\\
\displaystyle
\overline V(x,y,w,z)=\sup_{\xi\in\Xi}\inf_{\alpha\in{\cal A}} J(x,y,w,z,\alpha,\xi[\alpha]),
\end{array}
\end{equation}
\noindent
where $\Gamma$ and $\Xi$ are the set of non-anticipating strategies for player $X$ and player $Y$, respectively.

The main goal of the paper is to derive two suitable problems for Hamilton-Jacobi-Isaacs (HJI) equations in such a way to characterize $\underline V$ and $\overline V$ as the unique viscosity solutions of those problems, respectively. As a consequence, we will also get the existence of an equilibrium (i.e. $\underline V=\overline V$) under the standard Isaacs condition. To achieve the main goal we perform several steps. 

The first step is to prove the continuity in the space variables of the value functions, that is, for every $(w,z)\in\{-1,1\}\times\{-1,1\}$ fixed, the continuity of $(x,y)\mapsto V(x,y,w,z)$ (here and further by $V$ we denote any of two value functions (\ref{eq:values_introduction}), regardless whether it is the lower or the upper one). Under the hypothesis of decoupled dynamics and another decoupling hypothesis on the running cost $\ell$, such continuity is proved using a suitable construction of non-anticipating strategies. Indeed, in our switching differential game, we need the existence of some non-anticipating strategies which make the players, when they are on a switching threshold, to be able to switch or not (i.e. to cross the threshold or not) in dependence on its convenience. At the same time, such non-anticipating strategies must not penalize the cost too much. In the simpler case of an optimal control problem (one player only) this can be achieved by the Soner's construction of the so-called constrained controls \cite{son}. Indeed, the switching problem with state-dependent switching thresholds (as our problem is) is strongly related to state-constraints as well as exit-time problems. However, for the differential game situation, the stricter requirement on the construction of a Soner-like control that must be non-anticipating (i.e. non-dependent on future behaviors of the trajectories and controls), is a fundamental issue. Such an issue was addressed in the recent work by Bagagiolo-Maggistro-Zoppello \cite{bagmagzop}. In that work, the authors studied (for the first time) an exit-time/exit-costs differential game in the framework of dynamic programming and viscosity solutions theory for Isaacs equations with boundary conditions in the viscosity sense. In particular, the fundamental issue above is largely treated and solved under the decoupling hypotheses and further controllability hypotheses. In the present work, we apply results derived in \cite{bagmagzop}, where in the motivation part, we concerned the thermostatic problem, which is posed and
studied here.

The second step is to rewrite our infinite horizon problem as four exit-time/exit-costs problems coupled to each other by the exit costs. More precisely, the state space $\mathbb{R}^n\times\mathbb{R}^m\ni(X,Y)$ is divided in four partially overlapped sectors where the switching variables $(W,Z)$ remain constant. In each sector, the problem is considered as an exit-time differential game with exit costs mutually exchanged with the other sectors: the value function $V$ is evaluated in the new sector, after the previous one is left (i.e. one or both switching variable are switched), see Figure \ref{fig:sectors}. This step is achieved by dynamic programming techniques, using the decoupled feature of the dynamics, a controllability hypothesis and the already proved continuity of $V$.

Based on the second step, the third one is to write a system of four
HJI equations for four sectors, respectively, and coupled by the boundary conditions (for every sector, the boundary datum is the unknown function in the other sectors). Using the results of \cite{bagmagzop} on exit-time differential games, the value function $V$ will turn out to be a viscosity solution of such a system, where the boundary conditions are interpreted in a suitable viscosity sense.

The last step is to prove that $V$ is actually the unique solution of the system of HJI equations above. This is achieved by a fixed point procedure, using the uniqueness results of \cite{bagmagzop}. There, under some hypotheses, it is proved that the value function $V$ for a differential game of exit-time/exit-costs type (when the exit costs are given, not as here, where the exit costs are represented by the unknown solution itself) is the unique solution of a Dirichlet problem for the corresponding HJI equation. In particular, the boundary conditions are interpreted in a suitable viscosity sense that takes account of the min-max feature of the problem and benefits of the decoupling and controllability of the dynamics.


{\it Motivations and literature.}

There are different situations that can be interpreted as differential games with dynamics affected by switching. Just think to a pursuit evasion game (see Shinar-Glizer-Turetsky \cite{pursuer,evader}) where the switching dynamics is either the one of pursuer or the one of the evader. We can also imagine a race between two cars
where the switching variable(s) may represent the position of an automatic gears or the diesel/electric regime of an hybrid car as in Dextreit-Kolmanovsky \cite{Electric}.
 
Besides the well known shallow lake problem that arises in ecological
economics (see e.g., Reddy-Schumacher-Engwerda \cite{Engwerda}) can be seen as a differential game with switching dynamics as well as the international pollution problem with evolving environmental costs. Such costs, for less developed countries, change according to their cumulative revenue, see Masoudi-Zaccour \cite{Zaccour2013}.\\
We point out that the above mentioned switching dynamics are also called hybrid
dynamics.  A recent study of hybrid
differential games can be found also in Gromov-Gromova \cite{grogro}, where the authors formulated  necessary optimality conditions for determining optimal strategies in both cooperative and non-cooperative cases.
A particular class of differential games with changing structure is also considered in Bonneuil-Boucekkine \cite{Boucekkine}, where the transition to renewable energy leads to the change of the system's dynamics, and in Kort-Wrzaczek \cite{Kort}, where the change of a monopolist firm's dynamics is due to the entrance in the market of a firm offering the same products.
The switching can occurs not only in the dynamics but also in the cost function as,  for example, in Fabra-García \cite{Fabra} where a dynamic competition is analysed and the market prices change.

There are still mathematical motivations that suggest the study of differential games with thermostatic dynamics; these that are similar to the ones for studying optimal control
problems with thermostatic dynamics (see Bagagiolo-Maggistro\cite{bagmag}, Ceragioli-De Persis-Frasca \cite{CPF} and Bagagiolo-Danieli\cite{bagdan}).
In \cite{bagmag}, optimal controls problems
with dynamics inside a network are considered. A delay thermostat is introduced to overcome the discontinuity's problem arising when passing from an arc to another one due to the different dynamics and running cost on each branch of the network.
In \cite{CPF}, the authors make a rigorous treatment of continuous-time average consensus dynamics with uniform quantization in communications. The consensus is reached by quantized measurements which are transmitted using a delay thermostat.
Similarly, in \cite{bagdan} they consider an optimal control problem which has several internal switching variables that evolve following some delayed thermostatic laws.
A zero-sum differential games involving
hybrid controls was also considered in Dharmatti-Ramaswamy \cite{Ramaswamy}. Here, the state of the system is changed discontinuously and
the associated lower and upper value functions are characterized as the unique viscosity solutions of the corresponding quasi-variational
inequalities. Moreover, they give an Isaacs like condition for the game to have a value.

Up to the knowledge of the authors, the present work is the first attempt to study a switching/hybrid differential game in the framework of dynamic programming and viscosity solutions of Hamilton-Jacobi-Isaacs equations, and especially in connection with hybrid delayed thermostatic laws and, more in general, state-dependent switching.

In conclusion, we refer the reader to Bardi-Capuzzo Dolcetta \cite{barcap} for a comprehensive
account to viscosity solutions theory and applications to optimal control problems
and differential games (for differential games see also Buckdahn-Cardaliaguet-Quincampoix \cite{buccarqui}). Moreover, other
studies on constrained trajectories and non-anticipating strategies as well as on
possible relations with optimal control problems and differential games can be
found in Koike \cite{koi}, Bardi-Koike-Soravia \cite{barkoisor}, Cardaliaguet-Quincampoix-Saint
Pierre \cite{carquisain}, Bettiol-Cardaliaguet-Quincampoix \cite{betcarqui} Bettiol-Bressan-Vinter \cite{betbreal,betbrevin}, Bettiol-Facchi \cite{betfac} and Frankowska-Marchini-Mazzola \cite{framarmaz}.

{\it Plan of the paper.}

This paper is organized as follows: In Section 2, we introduce the hybrid thermostatic delayed relay, and its connection with ordinary differential equations. In Section 3, we briefly review the results in Bagagiolo-Maggistro-Zoppello \cite{bagmagzop} about exit-time/exit-costs differential games. In Section 4, we state the main assumptions about the infinite horizon switching differential game under study, and we argue about cost estimates on the switching trajectories. In section 5, we prove the continuity of the value functions and provide dynamic programming-like results, connecting the infinite horizon switching problem with exit-time problems. In Section 6, we characterize the value functions as the unique viscosity solutions of the systems of Isaacs equations. In Section 7, we give some hints on the numerical treatment of the problem providing a possible idea for a space-time discretization scheme. Section 8 draws conclusions and suggests future works.

\section{The hybrid thermostatic delayed relay}
\label{sec:relay}

\begin{figure}[htbp]
	\centering
	\includegraphics[scale=0.38]{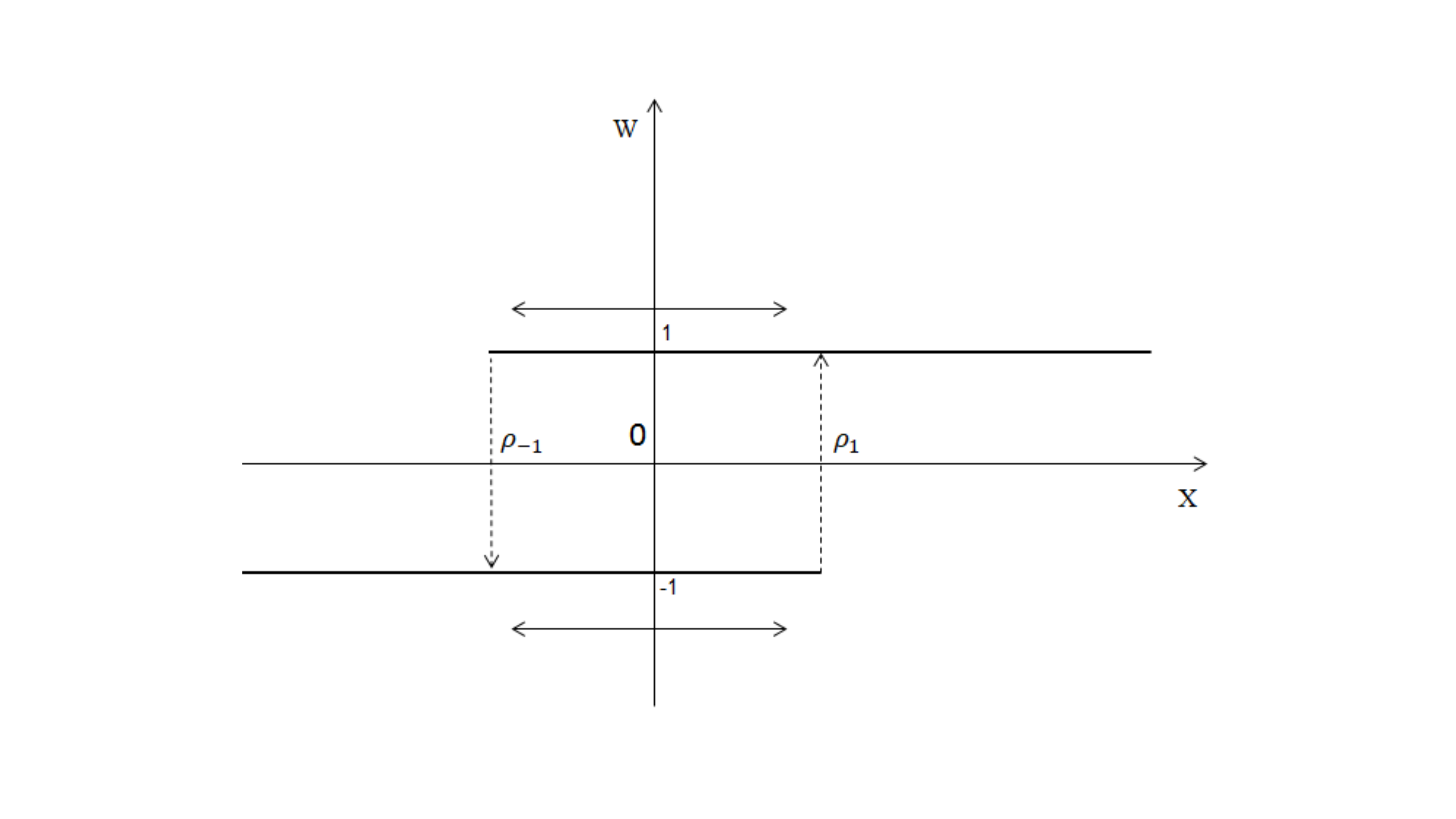}
	\vspace{-1cm}
	\caption{\label{fig:relay}Delayed thermostat with thresholds $\rho = (\rho_{-1}, \rho_1)$.}
\end{figure}

A hybrid delayed thermostat with thresholds $\rho=(\rho_{-1},\rho_1)$, $\rho_{-1}<\rho_1$, and initial output $w\in\{-1,1\}$, is the operator
\[
h_\rho[\cdot;w]:C^0(0,+\infty)\to L^\infty(0,+\infty),\ \ X\mapsto h_\rho[X;w],
\]
\noindent
whose behavior is described by Figure \ref{fig:relay}. In particular, it maps a time-continuous scalar input $X$ to a measurable time-dependent output function $W=h_\rho[X;w]$, which can only takes values in $\{-1,1\}$ and whose switching law is the following (see Visintin \cite{Visintin94} for a more systematic treatment of the delayed relay):
\[
\left\{
\begin{array}{ll}
\displaystyle
t\ge0,\ X(t)>\rho_1\Longrightarrow h_\rho[X;w](t)=1,\\
\displaystyle
t\ge0,\ X(t)<\rho_{-1}\Longrightarrow h_\rho[X;w](t)=-1,\\
\displaystyle
\rho_{-1}\le X(0)\le\rho_1\Longrightarrow\ h_\rho[X;w](0)=w,\\
\displaystyle
t\ge0,\ \rho_{-1}\le X(t)\le\rho_1\Longrightarrow\ h_{\rho}[X;w](t)=\lim_{s\to\tau_t^-}h_\rho[X;w](s)\\
\displaystyle
\ \ \ \ \mbox{if } \tau_t=\sup\{0\le\tau< t|X(\tau)<\rho_{-1}\ \mbox{or }\ X(\tau)>\rho_1\}>0,\\
\displaystyle
\ \ \ \ \ \ \ \ h_\rho[X;w](t)=h_\rho[X;w](0)=w\ \mbox{otherwise}.\\
\end{array}
\right.
\]
\noindent
In other words, looking to Figure \ref{fig:relay}, if at certain time $\tau$ we have $h_\rho[X;w](\tau)=-1$ (which certainly means $X(\tau)\le\rho_1$), then a switch from $-1$ to $1$ can only occur (and must occur) at a possible subsequent time $t\ge\tau$ if and only if, at that time $t$, the input $X$ crosses, strictly increasing, the upper threshold $\rho_1$, being $X(t)=\rho_1$. In that case, it is $h_\rho[X;w](t)=-1$ and $h_\rho[X;w]\equiv 1$ in, at least, a left-open right neighborhood of $t$, $]t,t+\delta]$. Similarly, if at certain time $\tau$ we have $h_\rho[X;w](\tau)=1$ (which certainly means $X(\tau)\ge\rho_{-1}$), then a switch from $1$ to $-1$ can only occur (and must occur) at a possible subsequent time $t\ge\tau$ if and only if, at that time $t$, the input $X$ crosses, strictly decreasing, the lower threshold $\rho_{-1}$, being $X(t)=\rho_{-1}$. In that case, it is $h_\rho[X;w](t)=1$ and $h_\rho[X;w]\equiv -1$ in, at least, a left-open right neighborhood of $t$, $]t,t+\delta]$. 
Note that, by such a definition, and, in particular, because of the strict inequality $\rho_{-1}<\rho_1$, the output $h_\rho[X;w]$ is left-continuous. In particular, we remark that, at a switching instant $t$, the value of the output $h_\rho[X;w](t)$ is still the previous one (i.e. it is not switched yet) and it will be equal to the new switched one at subsequent instants after $t$ only (if $X$ has crossed the threshold).
Finally note that the given initial output $w$ plays a role only if $\rho_{-1}\le X(0)\le\rho_1$.
An important property of the thermostatic delayed relay is the following semigroup property. For every $X$ and $w$, and for every $t,\tau\ge0$, it is

\begin{equation}
\label{eq:thermostat_semigroup}
h_\rho[X;w](t+\tau)=h_\rho[X(\cdot+t);h_\rho[X;w](t)](\tau).
\end{equation}

The switching evolution of $W=h_\rho[X;w]$ can be also described in the following way. Looking to Figure \ref{fig:relay}, we define

\[
{\cal O}_{-1}=]-\infty,\rho_1]\times\{-1\};\ \ {\cal O}_1=[\rho_{-1},+\infty[\times\{1\}
\]

\noindent
which correspond to the sets where the pair $(X,W)$ can evolve without switching. Given $X\in C^0(0,+\infty)$ and $w\in\{-1,1\}$ such that $(X(0),w)\in{\cal O}_{-1}\cup{\cal O}_1$, the output $h_\rho[X;w]$ is characterized as the unique left-continuous function $W$ such that

\[
\begin{array}{ll}
\displaystyle
(X(t),w(t))\in{\cal O}_{-1}\cup{\cal O}_1\ \forall\ t\ge0,\\
\displaystyle
\mbox{Var}_{[0,t]}W=\min\{\mbox{Var}_{[0,t]}\tilde W|(X(\tau),\tilde W(\tau))\in{\cal O}_{-1}\cup{\cal O}_1\forall\tau\in[0,t]\}\ \ \forall\ t\ge0,
\end{array}
\]

\noindent
where $\text{Var}_{[0,t]}$ is the total variation in the interval $[0,t]$ (which, for the delayed thermostat, corresponds to twice the number of switchings in $[0,t]$ (every switching has variation equal to $2$)). Hence, the switching law $t\mapsto W(t)$ is the unique one that satisfies the constraint $(X(t),W(t))\in{\cal O}_{-1}\cup{\cal O}_1$ and, in any time interval, minimizes the number of switchings. 

Such interpretation is also useful to define what is a solution of the switching scalar ordinary differential equation

\begin{equation}
\label{eq:thermostat_system}
\left\{
\begin{array}{ll}
\displaystyle
X'(t)=f(t,X(t),W(t)),\\
\displaystyle
W(t)=h_\rho[X;w](t),\\
\displaystyle
X(0)=x,\ \ (x,w)\in{\cal O}_{-1}\cup{\cal O}_1,
\end{array}
\right.
\end{equation}

\noindent
where $f:[0,+\infty[\times\mathbb{R}\times\{-1,1\}\to\mathbb{R},\ (t,x,w)\mapsto f(t,x,w)$ is bounded, measurable in $t\in[0,+\infty[$ and Lipschitz continuous in $x\in\mathbb{R}$ uniformly with respect to $(t,w)\in[0,+\infty[\times\{-1,1\}$. The solution is the unique function $t\mapsto(X(t),W(t))$ such that: i) $X$ is continuous and $W$ is left continuous; ii) $X(t)=x+\int_0^tf(s,X(s),W(s))ds$ and $(X(t),W(t))\in{\cal O}_{-1}\cup{\cal O}_1$ for all $t\ge0$; iii) minimizes in $[0,t]$ for all $t$, the number of switchings of $W$, among all pairs $(\tilde X,\tilde W)$ satisfying i) and ii). Such a unique solution can be constructed in the following way: consider the evolution of $X_1$, starting from $x$ with dynamics $f(\cdot,\cdot,w)$, and maintain such evolution until the possible time $t_1$ of switching for $W=h_\rho[X_1;w]$. Then in $]t_1,+\infty[$ consider the trajectory $X_2$ starting from $X_1(t_1)$ with dynamics $f(\cdot,\cdot,-w)$, and maintain it until the possible time $t_2>t_1$ of switching for $h_\rho[X_2,-w]$.  Then, in $]t_2,+\infty[$ consider the trajectory $X_3$ with dynamics $f(\cdot,\cdot,w)$ and starting from $X_2(t_2)$. Since the switching thresholds are different, $\rho_{-1}<\rho_1$, and the dynamics is bounded, then the number of those changes of dynamics is bounded in any compact set $[0,T]$ (there is not the so-called Zeno phenomenon: infinitely many switches in infinitesimal time intervals), and hence, gluing together the pieces of trajectories $X_i$, we get the unique solution of (\ref{eq:thermostat_system}) defined above. 

When the evolution is in $\mathbb{R}^n$ (as the controlled evolution of the next sections) we are going to suppose that the switching dependence of the dynamics is subject to the evolution of a fixed component of the continuous evolution $X\in\mathbb{R}^n$: $X\cdot\zeta$ where $\zeta\in\mathbb{R}^n$ is a unit vector. We then consider the $(n+1)$-dimensional system in the variable $(X,W)\in\mathbb{R}^n\times\{-1,1\}$

\[
\left\{
\begin{array}{ll}
\displaystyle
X'(t)=f(t,X(t),W(t)),\\
\displaystyle
W(t)=h_\rho[X\cdot\zeta;w](t),\\
\displaystyle
X(0)=x,\ \ (x\cdot\zeta,w)\in{\cal O}_{-1}\cup{\cal O}_1,
\end{array}
\right.
\]

\noindent
for which the same considerations as above hold.

\begin{Remark}
The previous definition of the switching rule $W=h_\rho[X;w]$, as well as the subsequent definition of the solution of (\ref{eq:thermostat_system}), is certainly linked to an exit-time feature: for example, a switch from $-1$ to $1$ occurs if and only if the pair $(X,W)$ exits from the closed set ${\cal O}_{-1}$. Indeed, in the following, we are going to interpret our infinite horizon switching differential games as a coupling of some exit-time differential games with exit from some suitable closed sets. One can also recast the problem as an exit-time problem with exit from an open set. Concerning the switching rule, it corresponds to the immediate switching when the threshold is touched (on the contrary, our definition is when the threshold is bypassed). This, for example,  would correspond to the fact that a switch from $-1$ to $1$ occurs if and only if the pair $(X,W)$ exits from the ``open" set $]-\infty,\rho_1[\times\{-1,1\}$. We use the interpretation as exit from a closed set because it is more prone to treat the Dirichlet problem from the Isaacs equation, since in such a case the boundary is ``physically" part of the problem, and it is even viable. However, assuming some suitable controllability conditions, as we are going to do, the two problems (exit from the open set and exit from the closed set) are in general extremely linked to each other, as the corresponding value functions in general coincide in the interior of the set. This fact is certainly known for optimal control problems but we do not investigate here such argument for the differential games situation.
\end{Remark}

\section{On exit-time/exit-costs differential games}
\label{sec:exit_time}

In this section, we briefly recall the results of Bagagiolo-Maggistro-Zoppello \cite{bagmagzop}, which will be used in the next sections.

Let us consider two open domains $\Omega_X\subseteq\mathbb{R}^n$, $\Omega_Y\subseteq\mathbb{R}^m$, with $C^2$ boundary and two decoupled controlled dynamics in $\mathbb{R}^n$ and $\mathbb{R}^m$, respectively

\begin{equation}
\label{eq:exit_system}
\left\{
\begin{array}{ll}
\displaystyle
X'(t)=f(X(t),\alpha(t))\\
\displaystyle
X(0)=x\in\overline\Omega_X
\end{array}
\right.,
\ \ \ \ \ \ 
\left\{
\begin{array}{ll}
\displaystyle
Y'(t)=g(Y(t),\beta(t))\\
\displaystyle
Y(0)=y\in\overline\Omega_Y
\end{array}
\right.
\end{equation}

\noindent
where for given compact sets $A,B$

\begin{equation}
\label{eq:measurable_controls}
\begin{array}{ll}
\displaystyle
\alpha\in{\cal A}=\{\alpha:[0,+\infty[\to A\ \mbox{measurable}\},\\
\displaystyle
\beta\in{\cal B}=\{\beta:[0,+\infty[\to B\ \mbox{measurable}\},\\
\displaystyle
f:\mathbb{R}^n\times A\to\mathbb{R}^n,\ \ \ g:\mathbb{R}^m\times B\to\mathbb{R}^m,\\
\displaystyle
f,g\ \mbox{are bounded and continuous and there exists } L>0\\
\displaystyle 
\mbox{such that for all } x_1,x_2\in\mathbb{R}^n,\ a\in A,\ y_1,y_2\in\mathbb{R}^m,\ b\in B\\
\displaystyle
\|f(x_1,a)-f(x_2,a)\|\le L\|x_1-x_2\|,\ \ \|g(y_1,b)-g(y_2,b)\|\le L\|y_1-y_2\|.
\end{array}
\end{equation}

We also consider the following functions

\begin{equation}
\label{eq:costs}
\begin{array}{ll}
\displaystyle
\ell:\mathbb{R}^n\times\mathbb{R}^m\times A\times B\to[0,+\infty[,\ (x,y,a,b)\mapsto\ell_1(x,y,a)+\ell_2(x,y,b),\\
\displaystyle
\ell_1:\mathbb{R}^n\times\mathbb{R}^m\times A\to[0,+\infty[,\ \ell_2:\mathbb{R}^n\times\mathbb{R}^m\times B\to[0,+\infty[,\\
\displaystyle
\Psi_X:\overline\Omega_X\to[0,+\infty[,\ \Psi_Y:\overline\Omega_Y\to[0,+\infty[,\ \Psi_{XY}:\partial\Omega_X\times\partial\Omega_Y\to[0,+\infty[,
\end{array}
\end{equation}

\noindent
with the assumption that $\ell_1,\ell_2,\Psi_X,\Psi_Y,\Psi_{XY}$ are bounded and continuous and that there exists $L>0$ such that for every $a\in A,b\in B$ fixed, and for all $x_1,x_2\in\mathbb{R}^n,\ y_1,y_2\in\mathbb{R}^m$,
\begin{equation}\label{lipelle}
\|\ell_1(x_1,y_1,a)-\ell_1(x_2,y_2,a)\|,\|\ell_2(x_1,y_1,b)-\ell_2(x_2,y_2,b)\|\le L\|(x_1,y_1)-(x_2,y_2)\|.
\end{equation}
We consider the differential game given by the cost, for $x\in\overline\Omega_X$, $y\in\overline\Omega_Y$,

\[
\begin{array}{ll}
\displaystyle
J(x,y,\alpha,\beta)=\int_0^{\tau_{(x,y)}(\alpha.\beta)}e^{-\lambda t}\ell(X(t),Y(t),\alpha(t),\beta(t))dt+\\
\displaystyle
e^{-\lambda\tau_{(x,y)}(\alpha,\beta)}\Psi(X(\tau_{(x,y)}(\alpha,\beta)),Y(\tau_{(x,y)}(\alpha,\beta))),
\end{array}
\]

\noindent
where $\lambda>0$ and, for $\alpha\in{\cal A},\beta\in{\cal B}$, $X(\cdot),Y(\cdot)$ are the trajectories given by (\ref{eq:exit_system}), sometimes also denoted as $X(\cdot;x,\alpha)$, $Y(\cdot;y,\beta)$,

\[
\begin{array}{ll}
\displaystyle
\tau_{(x,y)}(\alpha,\beta)=\min\{\tau_x(\alpha),\tau_y(\beta)\},\\
\displaystyle
\tau_x(\alpha)=\inf\{t\ge0\ \mbox{such that } X(t)\not\in\overline\Omega_X\},\\
\displaystyle
\tau_y(\beta)=\inf\{t\ge0\ \mbox{such that } Y(t)\not\in\overline\Omega_Y\},
\end{array}
\]

\noindent
and

\begin{equation}\label{psicosti}
\begin{array}{ll}
\displaystyle
\Psi(X(\tau_{(x,y)}(\alpha,\beta)),Y(\tau_{(x,y)}(\alpha,\beta)))=\\
\displaystyle
\left\{
\begin{array}{ll}
\displaystyle
\Psi_X(X(\tau_x(\alpha)),Y(\tau_x(\alpha)))&\mbox{if } \tau_{(x,y)}(\alpha,\beta)=\tau_x(\alpha)<\tau_y(\beta),\\
\displaystyle
\Psi_Y(X(\tau_y(\beta)),Y(\tau_y(\beta)))&\mbox{if } \tau_{(x,y)}(\alpha,\beta)=\tau_y(\beta)<\tau_x(\alpha),\\
\displaystyle
\Psi_{XY}(X(\tau_{(x,y)}(\alpha,\beta)),Y(\tau_{(x,y)}(\alpha,\beta)))&\mbox{if } \tau_{(x,y)}(\alpha,\beta)=\tau_x(\alpha)=\tau_y(\beta),\\
\end{array}
\right.
\end{array}
\end{equation}

\noindent
with, of course, $\inf\emptyset=+\infty$ and $e^{-\infty}\Psi:=0$.
Roughly speaking, the cost is paid as the integral of the discounted running cost up to the first exit-time of one of the two trajectories $X$ and $Y$ from its set of reference, $\overline\Omega_X$ and $\overline\Omega_Y$, respectively. Then a discounted exit cost is paid, which is given by three different exit costs, $\Psi_X,\Psi_Y,\Psi_{XY}$, depending whether player $X$ only exits from $\overline\Omega_X$ (i.e. $\tau_x<\tau_y$),  or player $Y$ only exits from $\overline\Omega_Y$ (i.e. $\tau_y<\tau_x$), or they both simultaneously exit from their closed reference sets (i.e. $\tau_x=\tau_y$).

The exit-time/exit-costs differential game is given by the fact that $X$ wants to minimize $J$ and $Y$ wants to maximize it. We then define the non-anticipating strategies for player $X$ and for player $Y$ respectively as:

\begin{equation}
\label{eq:strategies}
\begin{array}{ll}
\displaystyle
\Gamma=\Big\{\gamma:{\cal B}\to{\cal A}, \beta\mapsto\gamma[\beta]\ \mbox{such that }\\
\displaystyle
\beta_1=\beta_2\ \mbox{a. e. in } [0,t]\Longrightarrow\gamma[\beta_1]=\gamma[\beta_2]\ \mbox{a. e. in } [0,t],\ \forall t\ge0\Big\};\\
\displaystyle
\Xi=\Big\{\xi:{\cal A}\to{\cal B}, \alpha\mapsto\xi[\alpha]\ \mbox{such that }\\
\displaystyle
\alpha_1=\alpha_2\ \mbox{a. e. in } [0,t]\Longrightarrow\xi[\alpha_1]=\xi[\alpha_2]\ \mbox{a. e. in } [0,t],\ \forall t\ge0\Big\}.
\end{array}
\end{equation}

The lower and upper value functions are respectively defined as, for $(x,y)\in\overline\Omega_X\times\overline\Omega_Y$, 

\begin{equation} \label{lowUpp}
\begin{array}{ll}
\displaystyle
\underline v(x,y)=\inf_{\gamma\in\Gamma}\sup_{\beta\in{\cal B}}J(x,y,\gamma[\beta],\beta),\\
\displaystyle
\overline v(x,y)=\sup_{\xi\in\Xi}\inf_{\alpha\in{\cal A}}J(x,y,\alpha,\xi[\alpha]).
\end{array}
\end{equation}

Similarly to the non-anticipating strategies, we define a non-anticipating tuning for both players: a non-anticipating tuning is any function ${\cal K}\to{\cal K}, k\mapsto\tilde k$, where ${\cal K}$ is either $\cal A$ or $\cal B$, such that

\begin{equation}
\label{eq:tuning}
k_1=k_2\ \mbox{a. e. in } [0,t]\ \Longrightarrow\ \tilde k_1=\tilde k_2\ \mbox{a. e. in } [0,t],\ \forall\ t\ge0.
\end{equation}

Note the difference: a non-anticipating strategy is a function from the set of measurable controls for one player to the set of measurable controls for the other player; a non-anticipating tuning is a function from the set of measurable controls for one player to itself.

We then assume the following controllability and compatibility hypotheses

\begin{equation}
\label{eq:assumptions}
\begin{array}{ll}
\displaystyle
\forall\ x\in\partial\Omega_X\ \exists\ a_1,a_2\in A\ \mbox{such that } f(x,a_1)\cdot\ n_X(x)<0<f(x,a_2)\cdot n_X(x),\\
\displaystyle
\forall\ y\in\partial\Omega_Y\ \exists\ b_1,b_2\in B\ \mbox{such that } g(y,b_1)\cdot\ n_Y(x)<0<g(y,b_2)\cdot n_Y(y),\\
\displaystyle
\Psi_Y(x,y)\le\Psi_{XY}(x,y)\le\Psi_X(x,y)\ \ \forall\ (x,y)\in\partial\Omega_X\times\partial\Omega_Y,
\end{array}
\end{equation}

\noindent
where $n_X(x)$ and $n_Y(y)$ are, respectively, the outer normal unit vector to $\overline\Omega_X$ in $x$ and to $\overline\Omega_Y$ in $y$.

Finally, we define, respectively, the upper Hamiltonian and the lower Hamiltonian for $(x,y,p,q)\in\mathbb{R}^n\times\mathbb{R}^m\times\mathbb{R}^n\times\mathbb{R}^m$, as

\[
\begin{array}{ll}
\displaystyle
UH(x,y,p,q)=\min_{b\in B}\max_{a\in A}\{-f(x,a)\cdot p-g(y,b)\cdot q-\ell(x,y,a,b)\},\\
\displaystyle
LH(x,y,p,q)=\max_{a\in A}\min_{b\in B}\{-f(x,a)\cdot p-g(y,b)\cdot q-\ell(x,y,a,b)\}.
\end{array}
\]

\begin{Theorem}
\label{thm:bagmagzop}
Under hypotheses \eqref{eq:measurable_controls}-\eqref{eq:assumptions}, the value functions $\underline v$ and $\overline v$ \eqref{lowUpp} are bounded and continuous in $\overline\Omega_X\times\overline\Omega_Y$. Moreover, they are, respectively, the unique bounded and continuous function $u:\overline\Omega_X\times\overline\Omega_Y\to\mathbb{R}$ which satisfies, in the viscosity sense, the following Dirichlet problems for the Isaacs equations

\begin{equation}
\label{eq:lower_dirichlet}
\left\{
\begin{array}{ll}
\displaystyle
\lambda u(x,y)+UH(x,y,\nabla_xu(x,y),\nabla_yu(x,u))=0&\mbox{in } \Omega_X\times\Omega_Y,\\
\displaystyle
u=\Psi_X&\mbox{on } \partial\Omega_X\times\Omega_Y\\
\displaystyle
u=\Psi_Y&\mbox{on } \Omega_X\times\partial\Omega_Y\\
\displaystyle
u=\Psi_X\ \mbox{or } u=\Psi_Y&\mbox{on } \partial\Omega_X\times\partial\Omega_Y
\end{array}
\right.
\end{equation}

\begin{equation}
\label{eq:upper_dirichlet}
\left\{
\begin{array}{ll}
\displaystyle
\lambda u(x,y)+LH(x,y,\nabla_xu(x,y),\nabla_y(x,y))=0&\mbox{in } \Omega_X\times\Omega_Y,\\
\displaystyle
u=\Psi_X&\mbox{on } \partial\Omega_X\times\Omega_Y\\
\displaystyle
u=\Psi_Y&\mbox{on } \Omega_X\times\partial\Omega_Y\\
\displaystyle
u=\Psi_X\ \mbox{or } u=\Psi_Y&\mbox{on } \partial\Omega_X\times\partial\Omega_Y
\end{array}
\right.
\end{equation}

\noindent
where $\nabla_x$ and $\nabla_y$ stay, respectively, for the gradient with respect to the $x\in\mathbb{R}^n$ variable, and the gradient with respect to the $y\in\mathbb{R}^m$ variable.
\end{Theorem}

In \cite{bagmagzop}, an ad-hoc definition of viscosity solution is given and used, especially for what concerns the boundary conditions, in order to suitably treat the min-max feature of the problem and the separation of the three exist costs. By a solution in the viscosity sense of the problem (\ref{eq:lower_dirichlet}) (and similarly for (\ref{eq:upper_dirichlet})), we mean the following: let $\varphi\in C^1(\overline\Omega_X\times\overline\Omega_Y)$ be a test function, and $(x_0,y_0)\in\overline\Omega_X\times\overline\Omega_Y$, then the following facts i) and ii) hold true:

i) if $(x_0,y_0)$ is a point of local maximum for $u-\varphi$, with respect to $\overline\Omega_X\times\overline\Omega_Y$, then we have the following four implications (one per every line)

\begin{equation}
\label{eq:subsol}
\begin{array}{l}
\left.
\begin{array}{ll}
\displaystyle
(x_0,y_0)\in\Omega_X\times\Omega_Y,\\
\displaystyle
(x_0,y_0)\in\partial\Omega_X\times\Omega_Y,\ u(x_0,y_0)>\psi_X(x_0,y_0),\\
\displaystyle
(x_0,y_0)\in\Omega_X\times\partial\Omega_Y,\ u(x_0,y_0)>\psi_Y(x_0,y_0),\\
\displaystyle
(x_0,y_0)\in\partial\Omega_X\times\partial\Omega_Y,\  \psi_X(x_0,y_0)\neq u(x_0,y_0)>\psi_Y(x_0,y_0)\\
\end{array}
\right\}\ \Longrightarrow\\
\\
\displaystyle
\ \ \ \ \lambda u(x_0,y_0)+UH(x_0,y_0,\varphi_x(x_0,y_0),\varphi_y(x_0,y_0))\le0;
\end{array}
\end{equation}

ii)  if $(x_0,y_0)$ is a point of local minimum for $u-\varphi$, with respect to $\overline\Omega_X\times\overline\Omega_Y$, then we have the following four implications (one per every line)

\begin{equation}
\label{eq:supersol}
\begin{array}{l}
\left.
\begin{array}{ll}
\displaystyle
(x_0,y_0)\in\Omega_X\times\Omega_Y,\\
\displaystyle
(x_0,y_0)\in\partial\Omega_X\times\Omega_Y,\ u(x_0,y_0)<\psi_X(x_0,y_0),\\
\displaystyle
(x_0,y_0)\in\Omega_X\times\partial\Omega_Y,\ u(x_0,y_0)<\psi_Y(x_0,y_0),\\
\displaystyle
(x_0,y_0)\in\partial\Omega_X\times\partial\Omega_Y,\  \psi_Y(x_0,y_0)\neq u(x_0,y_0)<\psi_X(x_0,y_0)\\
\end{array}
\right\}\ \Longrightarrow\\
\\
\displaystyle
\ \ \ \ \lambda u(x_0,y_0)+UH(x_0,y_0,\varphi_x(x_0,y_0),\varphi_y(x_0,y_0))\ge0;
\end{array}
\end{equation}

The implications given by the second, third and fourth lines of (\ref{eq:subsol})--(\ref{eq:supersol}) represent the boundary conditions in the viscosity sense.

\begin{Remark}
\label{rmrk:tuning}
Note that in the definition of the boundary conditions in viscosity sense above, the exit cost $\Psi_{XY}$ for the simultaneous exit of both players, actually does not play any role. This is a consequence of the compatibility condition in \eqref{eq:assumptions}, which is, in some sense, a sort of stability: the exit cost for the minimizing player is larger than the cost of the maximizing one.

The decoupled feature of the dynamics and of the running cost, together with the controllability and compatibility conditions in (12) as well as the regularity of the boundary, play and important role for the continuity and uniqueness result of Theorem 2.

In particular, for what concerns the continuity, in \cite{bagmagzop} it is proved that, under hypotheses \eqref{eq:measurable_controls}-\eqref{eq:assumptions}, the following property holds:

for every $T>0$, for every $K\subset\mathbb{R}^n\times\mathbb{R}^m$ compact, there exist $\delta>0$ and a modulus of continuity $\mathcal{O}_{T,K}$, and: 

I) for every $(x_1,y_1),(x_2,y_2)\in K\cap(\overline\Omega_X\times\overline\Omega_Y)$, with $\|(x_1,y_1)-(x_2,y_2)\|\le\delta$ there exists a non-anticipating tuning $\beta\mapsto\overline\beta$ from $\cal B$ to itself, and there exists a way to associate $\overline\gamma\in\Gamma$ to any $\gamma\in\Gamma$, such that, for every $\beta\in{\cal B}$, $\gamma\in\Gamma$, we have

\begin{equation}
\label{eq:Assumption2}
\begin{array}{l}
\displaystyle
i)\ 0\le\tau_{x_1}(\overline\gamma[\beta])-\tau_{x_2}(\gamma[\beta])\le{\cal O}_{T,K}(\|x_1-x_2\|),\\
\displaystyle
ii)\ 0\le\tau_{y_2}(\overline\beta)-\tau_{y_1}(\beta)\le {\cal O}_{T,K}(\|y_1-y_2\|),\\
\displaystyle
iii)\ \|X(t;x_1,\overline\gamma[\beta])-X(t;x_2,\gamma[\beta])\|\le\mathcal{O}_{T,K}(\|x_1-x_2\|),\ \forall\ t\in[0,\tilde\tau]\\
\displaystyle
iv)\ \|Y(t;y_1,\beta)-Y(t;y_2,\overline\beta)\|\le\mathcal{O}_{T,K}(\|y_1-y_2\|),\ \forall\ t\in[0,\tilde\tau]\\
\displaystyle
v)\ \left|J_{\tilde\tau}(x_1,y_1,\overline\gamma[\overline\beta],\beta)-J_{\tilde\tau}(x_2,y_2,\gamma[\overline\beta],\overline\beta)\right|\le\mathcal{O}_{T,K}(\|(x_1,y_1)-(x_2,y_2)\|),\\
\end{array}
\end{equation}

\noindent
where $\tilde\tau=\min\{\tau_{x_2}(\gamma[\overline\beta]),\tau_{y_1}(\beta),T\}$, and $J_{\tilde\tau}$ is the integral of the discounted running cost up to the time $\tilde\tau$: $\int_0^{\tilde\tau}e^{-\lambda t}\ell dt$.

II) Similarly it holds reversing the roles of $X$ and $Y$, $\gamma\in\Gamma$ and $\xi\in\Xi$, $\alpha\in\cal A$ and $\beta\in\cal B$.

We point out that the construction of $\overline\beta$ and of $\overline\gamma$ are made independently on the behavior of the other player. That is, $\overline\beta(t)$ and $\overline\gamma[\beta]$ are dependent only on the behavior, up to the time $t$, of the trajectories $Y(\cdot;y_1,\beta)$ and $X(\cdot;x_2,\gamma[\beta])$, respectively. This is possible essentially due to the decoupling feature in the controls of the running cost $\ell$ \eqref{eq:costs} (see \cite{bagmagzop}, section 7, Remark 12).

Note that in \eqref{eq:Assumption2}, i) means that the trajectory starting from $x_1$ with control $\overline\gamma[\beta]$ does not exit before the trajectory starting from $x_2$ with control $\gamma[\beta]$, and moreover, the difference of the two exit instants are controlled by the initial distance of the points; ii) means that the trajectory starting from $y_2$ with control $\overline\beta$ does not exit before the trajectory starting from $y_1$ with control $\beta$, and the difference is controlled by the initial distance; iii), iv), v) mean that the distance of those (and other similar) trajectories and their costs are controlled by the initial distances. Of course, if $\tau_{x_2}(\gamma[\beta])=+\infty$, then both trajectories never exit, and similarly if $\tau_{y_1}(\beta)=+\infty$.

This property is essential in order to prove the continuity. Under the hypotheses here stated, its validity is proven in \cite{bagmagzop} (see Assumption 2, points 3) and 7) of Proposition 3, and (7)--(10)), suitably adapting the construction in Soner \cite{son} to the non-anticipating framework.

We are going to use \eqref{eq:Assumption2} in the next sections. 
\end{Remark}

\section{The switching infinite horizon differential game}
\label{sec:assumptions}
The decoupled controlled dynamics of the players are respectively given by

\begin{equation}
\small
\label{eq:thermostat_systems}
\left\{
\begin{array}{ll}
X'(t)=f(X(t),W(t),\alpha(t)),\ t>0\\
W(t)=h_\rho[X\cdot\zeta_X;w](t)\\
(X(0),w(0))=(x,w)
\end{array}
\right.,
\ \ \ 
\left\{
\begin{array}{ll}
Y'(t)=g(Y(t),Z(t),\beta(t)),\ t>0\\
Z(t)=h_\eta[Y\cdot\zeta_Y;z](t)\\
(Y(0),z(0))=(y,z)
\end{array}
\right.
\end{equation}
\noindent
where each dynamics is affected by a delayed thermostatic switching rule. 

Here and in the sequel we will assume the following hypotheses and use the following notations:\\

{\it Main Assumptions}
\begin{itemize}
\item $X(t)=X(t;x,w,\alpha)\in\mathbb{R}^n$, $Y(t)=Y(t;y,z,\beta)\in\mathbb{R}^m$ are the states at time $t$ of the player $X$ and player $Y$ whose evolution are given by the trajectories of (\ref{eq:thermostat_systems}), respectively (here and in the sequel, the names of the players will be identified with the names of their state variable);
\item $h_\rho$ and $h_\eta$ are delayed switching thermostats with thresholds $\rho_{-1}<\rho_1$ and $\eta_{-1}<\eta_1$, respectively;
\item $W(t)=W(t;x,w,\alpha)\in\{-1,1\}$ and $Z(t)=Z(t;y,z,\beta)\in\{-1,1\}$ are the switching variables, with evolutions given by (\ref{eq:thermostat_systems}), respectively; 
\item $\zeta_X\in\mathbb{R}^n$, $\zeta_Y\in\mathbb{R}^m$ are unit vectors;  $X\cdot\zeta_X$, $Y\cdot\zeta_Y$ are scalar products in $\mathbb{R}^n$ and $\mathbb{R}^m$, respectively, and represent the input functions $t\mapsto X(t)\cdot\zeta_X$, $t\mapsto Y(t)\cdot\zeta_Y$ to which, via the delayed thermostats, the switching laws of the variables $W$ and $Z$ are subject;
\item $A$, $B$ (the sets of constants controls), $\cal A$, $\cal B$ (the sets of measurable controls), $\Gamma$ and $\Xi$ (the sets of non-anticipating strategies) are defined as in (\ref{eq:measurable_controls}) and (\ref{eq:strategies}); a non-anticipating tuning $k\mapsto\tilde k$ is defined as in (\ref{eq:tuning});
\item $(x,w)\in\mathbb{R}^n\times\{-1,1\}$, $(y,z)\in\mathbb{R}^m\times\{-1,1\}$ are suitable initial data;
\item $f:\mathbb{R}^n\times\{-1,1\}\times A\to\mathbb{R}^n$, $g:\mathbb{R}^m\times\{-1,1\}\times B\to\mathbb{R}^m$ are the switching controlled dynamics of player $X$ and player $Y$, respectively. Moreover, they are continuous, bounded and Lipschitz in the state variables, i.e. $$\exists\ M>0\ \mbox{such that}\ \forall (x,y,w,z,a,b)\\ 
\displaystyle
\quad\|f(x,w,a)\|,\|g(y,z,b)\|\le M,$$ and 
$$
\begin{aligned}
&\exists\ L>0\ \mbox{such that}\ \forall (x_1,w,a),(x_2,w,a),(y_1,z,b),(y_2,z,b)\\
\displaystyle
&\|f(x_1,w,a)-f(x_2,w,a)\|\le L\|x_1-x_2\|,\\
\displaystyle
&|g(y_1,z,b)-g(y_2,z,b)\|\le L\|y_1-y_2\|.
\end{aligned}
$$
\item $\ell:\mathbb{R}^n\times\mathbb{R}^m\times\{-1,1\}\times\{-1,1\}\times A\times B\to[0,+\infty[$, $(x,y,w,z,a,b)\mapsto\ell(x,y,w,z,a,b)=\ell_1(x,y,w,z,a)+\ell_2(x,y,w,z,b)$ is the running cost, decoupled in the controls, where $\ell_1$ and $\ell_2$ are continuous, bounded and Lipschitz continuous with respect to the state variables; in particular $\exists\ M,L>0$ such that $\forall x,y,x_1,y_1,x_2,y_2,w,z,a,b$:

\[
\begin{array}{ll}
\displaystyle
\|\ell_1(x,y,w,z,a)\|,\ \|\ell_2(x,y,w,z,b)\|\le M,\\
\displaystyle
\|\ell_1(x_1,y_1,w,z,a)-\ell_1(x_2,y_2,w,z,a)\|\le L\|(x_1,y_1)-(x_2,y_2)\|,\\
\displaystyle
\|\ell_2(x_1,y_1,w,z,b)-\ell_2(x_2,y_2,w,z,b)\|\le L\|(x_1,y_1)-(x_2,y_2)\|.
\end{array}
\]

\item $\lambda>0$ is the discount factor.
\end{itemize}

We consider an infinite horizon discounted problem where, as usual, $X$ wants to minimize and $Y$ wants to maximize a cost of the form
\[
J(x,y,w,z,\alpha,\beta)=\int_0^{+\infty}e^{-\lambda t}\ell(X(t),Y(t),W(t),Z(t),\alpha(t),\beta(t))dt.
\]

\noindent
Note that the cost $J$ is also depending on the switching variables $W$ and $Z$.

\begin{figure}[H]
	\label{fig:sectors}
	\centering
	\includegraphics[scale=0.47]{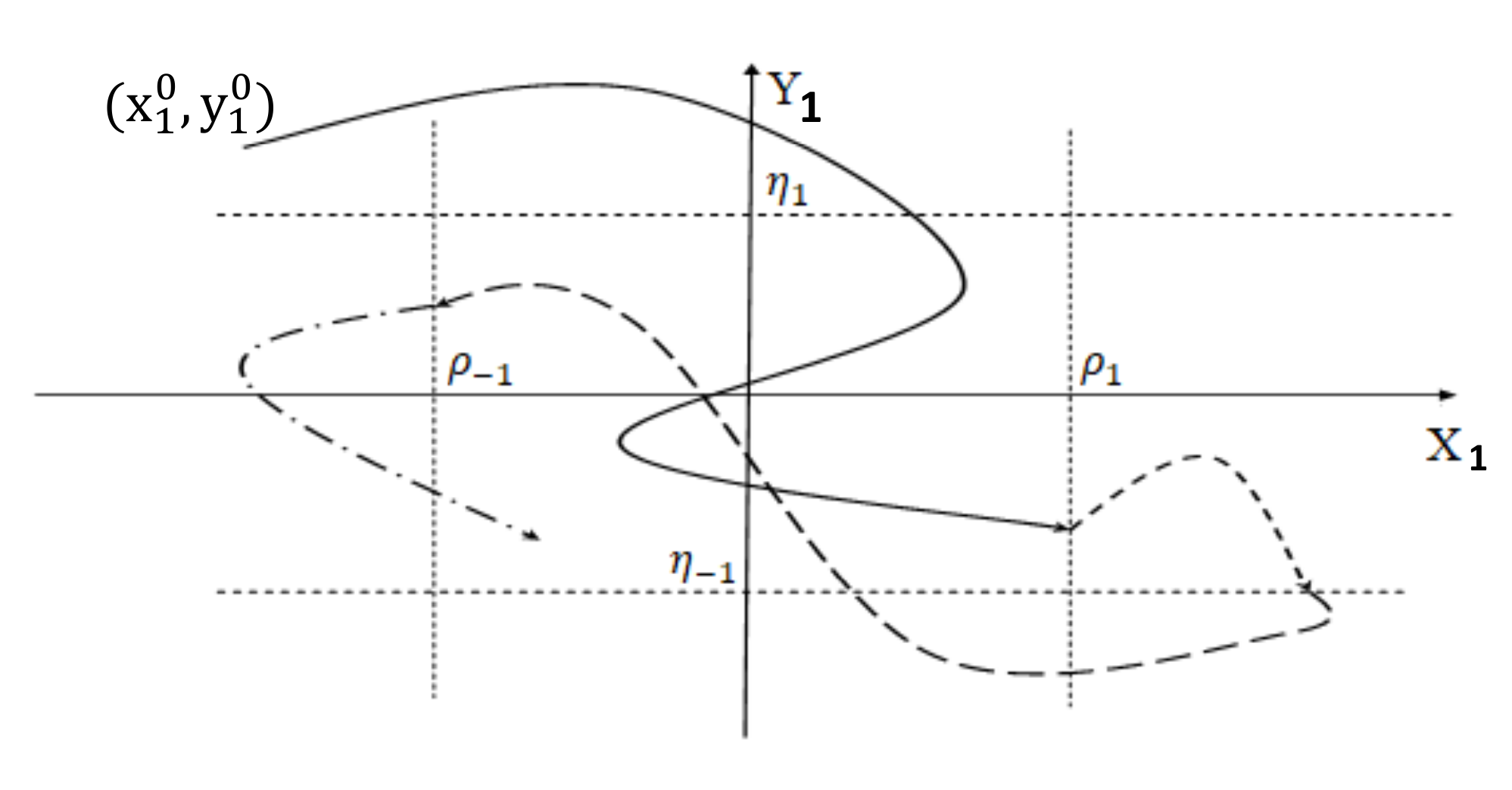}
	\caption{\label{projection}The trajectory of the pair $(X_1,Y_1)$ in the plane}
\end{figure}

We then define the lower and upper value functions as, respectively
\begin{equation}
\label{eq:values}
\begin{array}{ll}
\displaystyle
\underline V(x,y,w,z)=\inf_{\gamma\in\Gamma}\sup_{\beta\in{\cal B}} J(x,y,w,z,\gamma[\beta],\beta),\\
\displaystyle
\overline V(x,y,w,z)=\sup_{\xi\in\Xi}\inf_{\alpha\in{\cal A}} J(x,y,w,z,\alpha,\xi[\alpha]).
\end{array}
\end{equation}
\noindent
In order to simplify notations, we assume that $\zeta_X$ and $\zeta_Y$ are the first unit canonical vectors, so that $X_1=X\cdot\zeta_X$ and $Y_1=Y\cdot\zeta_Y$ are the first coordinates of $X$ and $Y$, respectively.\\
Let us consider the evolution of $X$ given by (\ref{eq:thermostat_systems}). We can interpret such an evolution as a switching evolution governed by two dynamics-modes, $f(\cdot,1,\cdot)$ and $f(\cdot,-1,\cdot)$, where the switching between the two modes is governed by the delayed thermostat $h_\rho$ subject to the evolution of $X_1$. Similarly, the evolution of the player $Y$ given by \eqref{eq:thermostat_systems}, which is affected by the delayed thermostat $h_\eta$ subject to the evolution of $Y_1$, switches between the two dynamics $g(\cdot, 1,\cdot)$ and $g(\cdot,-1,\cdot)$. The behavior of the projection on the first coordinates of $X$ and $Y$ respectively,  is described by Figure \ref{projection}. For example, for given controls $\alpha\in{\cal A}$, and $\beta\in{\cal B}$,  the filled curve is the evolution with dynamics $(f(\cdot,-1,\alpha), g(\cdot,1,\beta))$, the short dashed curve is the evolution with $(f(\cdot,1,\alpha), g(\cdot,1,\beta))$, the long dashed curve is the evolution with $(f(\cdot,1,\alpha), g(\cdot,-1,\beta))$ and the point-dashed one is the evolution with $(f(\cdot,-1,\alpha), g(\cdot,-1,\beta))$.
As it is easily deduced by the trajectories described in Figure \ref{projection}, the state space $\mathbb{R}^n\times\mathbb{R}^m$  can be divided in $4$ (non-disjointed, but overlapped) closed sectors, every one indexed by the corresponding $2$-string $(w,z)$ of $1$ and $-1$. More precisely
$$
\mathbb{R}^n\times\mathbb{R}^m=\bigcup_{(w,z)\in\{-1,1\}^2}\overline{B}_{(w,z)}=\Bigl(\overline{B}_{(1,1)}\cup\overline{B}_{(1,-1)}\cup\overline{B}_{(-1,1)}\cup\overline{B}_{(-1,-1)}\Bigr),
$$
where
\begin{equation}
\label{eq:B}
\begin{aligned}
&\overline{B}_{(1,1)}=[\rho_{-1},+\infty[\times\mathbb{R}^{n-1}\times[\eta_{-1},+\infty[\times\mathbb{R}^{m-1},\\
&\overline{B}_{(1,-1)}=[\rho_{-1},+\infty[\times\mathbb{R}^{n-1}\times]-\infty,\eta_1]\times\mathbb{R}^{m-1},\\
&\overline{B}_{(-1,1)}=]-\infty,\rho_1]\times\mathbb{R}^{n-1}\times[\eta_{-1},+\infty[\times\mathbb{R}^{m-1},\\
&\overline{B}_{(-1,-1)}=]-\infty,\rho_1]\times\mathbb{R}^{n-1}\times]-\infty,\eta_{1}]\times\mathbb{R}^{m-1}.
\end{aligned}
\end{equation}
When we start to move inside one of the sectors, then we continue to move in the same mode $(f(\cdot,w,\cdot),g(\cdot,z,\cdot))$ until we leave that sector, and after that we move in the new modality (corresponding to the index of the new sector) determined by the delayed thermostatic switching rules $h_\rho$ and $h_\eta$.

In the next section, we are going to interpret the switching infinite horizon problem as four exit-time/exit-costs problems, one per every sector, and coupled by mutually exchanged exit-costs.  In order to recast such exit-time/exit-costs problems in the framework of Section \ref{sec:exit_time}, let us note that every sector is of the form $\overline B_{(w,z)}=\overline\Omega_X^w\times\overline\Omega_Y^z$ where, for example,
$\overline\Omega_X^1=[\rho_{-1},+\infty[\times\mathbb{R}^{n-1}$ and $\overline\Omega_Y^1=[\eta_{-1},+\infty[\times\mathbb{R}^{m-1}$. Note that $\partial\Omega_X^w=\{\rho_{-w}\}\times\mathbb{R}^{n-1}$ and $\partial\Omega_Y^z=\{\eta_{-z}\}\times\mathbb{R}^{m-1}$, that is the boundaries are the switching thresholds points (and the threshold is crossed when the first component crosses it). Moreover, we also require the following controllability assumption

\begin{equation}
\label{eq:controllability-switching}
\begin{array}{ll}
\displaystyle
\forall\ w\in\{-1,1\},\ \forall x\in\partial\Omega_X^w\ \exists\ a_1,a_2\in A\\
\displaystyle
\ \ \ \ \mbox{such that } f_1(x,w,a_1)<0<f_1(x,w,a_2),\\
\displaystyle
\forall\ z\in\{-1,1\},\ \forall y\in\partial\Omega_Y^z\ \exists\ b_1,b_2\in B\\
\displaystyle
\ \ \ \ \mbox{such that } g_1(y,z,b_1)<0<g_1(y,z,b_2),
\end{array}
\end{equation}

\noindent
where $f_1$ and $g_1$ are, respectively, the first components of the dynamics $f$ and $g$. Note that (\ref{eq:controllability-switching}) means that, when $X$ or $Y$ are in a switching threshold, then, due to the decoupled feature of the dynamics, they can freely choose whether to switch or not (remember that a switch occurs only when the threshold is bypassed). 



Finally, note that whenever a finite time $T > 0$ is fixed, then any pair of switching trajectories $(W,Z)$ given by (\ref{eq:thermostat_systems}) can switch only an a-priori bounded finite number of times in $[0,T]$. This is true because the dynamics are bounded, and the thresholds are disjointed: for example, the trajectory $X_1$ needs a uniform positive time $t>0$ in order to pass from $\rho_{-1}$ to $\rho_1$ and vice-versa. Let $N_T>0$ be such an a-priori bound for the number of switches in $[0,T]$. For every $(x,y,w,z)$ and every control $\alpha\in{\cal A},\beta\in{\cal B}$, and for every $T>0$ we then have a finite sequence of switching instants (possibly empty, if the trajectories never switch) in $[0,T]$:

\begin{equation}
\label{eq:sequences}
\begin{array}{ll}
\displaystyle
\tau^1_X(x,w,\alpha)<\tau^2_X(x,w,\alpha)<\cdots<\tau^{N_X}_X(x,w,\alpha),\\
\displaystyle
\tau^1_Y(y,z,\beta)<\tau^2_Y(y,z,\beta)<\cdots<\tau^{N_Y}_Y(y,z,\beta),\\
\displaystyle
\tau^1_{XY}(x,y,w,z,\alpha,\beta)<\tau^2_{XY}(x,y,w,z,\alpha,\beta)<\cdots<\tau^{N_{XY}}_{XY}(x,y,w,z,\alpha,\beta),
\end{array}
\end{equation}

\noindent
with $N=N_X+N_Y+N_{XY}\le N_T$, and where $\tau_X$ corresponds to switches of $X$ only, $\tau_Y$ to switches of $Y$ only, and $\tau_{XY}$ to simultaneous switches of $X$ and $Y$ (which means $\tau_X=\tau_Y$). We can merge such three sequences, in order to get a unique sequence

\begin{equation}
\label{eq:unique_sequence}
\tau^0:=0\le\tau^1<\tau^2<\cdots<\tau^N\le T=:\tau^{N+1},\ \ N\le N_T.
\end{equation}

\noindent
Moreover, we denote by $(w,z)^i$ the new values of the switching variables after the $i$-th switch, $i=1,\dots,N$, and define $(w,z)^0=(w,z)$. We have

\begin{equation}
\label{eq:sum_costs}
\begin{array}{ll}
\displaystyle
J_T(x,y,w,z,\alpha,\beta)=\\
\displaystyle
\sum_{i=1}^{N+1}\int_{\tau^{i-1}}^{\tau^i}e^{-\lambda t}\ell\left(X(t),Y(t),(w,z)^{i-1},\alpha(t),\beta(t)\right)dt=\\
\displaystyle
\sum_{i=1}^{N+1}e^{-\lambda\tau^{i-1}}\int_0^{\tau^i-\tau^{i-1}}e^{-\lambda t}\ell\left(X^{i-1}(t),Y^{i-1}(t),(w,z)^{i-1},\alpha^{i-1},\beta^{i-1}\right)dt=\\
\displaystyle
\sum_{i=1}^{N+1}e^{-\lambda\tau^{i-1}}J_{(\tau^i-\tau^{i-1})}\left(X(\tau^{i-1}),Y(\tau^{i-1}),(w,z)^{i-1},\alpha^{i-1},\beta^{i-1}\right),
\end{array}
\end{equation}

\noindent
where $(X^{i-1},Y^{i-1})$ is the trajectory starting from $(X(\tau^{i-1}),Y(\tau^{i-1}),(w,z)^{i-1})$ with controls $(\alpha^{i-1}(\cdot),\beta^{i-1}(\cdot))=(\alpha(\cdot+\tau^{i-1}),\beta(\cdot+\tau^{i-1}))$.

Using the representation (\ref{eq:sum_costs}), in the spirit of (\ref{eq:Assumption2}), we now construct a non-anticipating tuning and a non-anticipating strategy, which work for our switching problem.
Take $\mu>0$ such that if a trajectory switches at time $\tau$, then it does not switch in the time interval $[\tau,\tau+\mu]$.
Take $T>0$ and $K\subset\mathbb{R}^n\times\mathbb{R}^m$ compact.
Fix $(w,z)$ and $(x_1,y_1),(x_2,y_2)\in K\cap\overline B_{(w,z)}$ such that $\|(x_1,y_1)-(x_2,y_2)\|\le\delta$, where $\delta$ is as in (\ref{eq:Assumption2}), and moreover, such that ${\cal O}_{T,K}(e^{LT}\delta)<\mu/2$. Take $\beta\in{\cal B}$ and $\gamma\in\Gamma$. By (\ref{eq:Assumption2}), with the notations of (\ref{eq:sequences}), we get the non-anticipating tuning $\overline\beta^0$ and the non-anticipating strategy $\overline\gamma^0$ such that

\begin{equation}
\label{eq:Assumption2-type}
\begin{array}{ll}
\displaystyle
0\le\tau^1_Y(y_2,z,\overline{\beta}^0)-\tau^1_Y(y_1,z,\beta)\le{\cal O}_{T,K}(\|y_1-y_2\|),\\
\displaystyle
0\le\tau^1_X(x_1,w,\overline{\gamma}^0[\beta])-\tau^1_X(x_2,w,\gamma[\beta])\le{\cal O}_{T,K}(\|x_1-x_2\|),\\
\displaystyle
\|Y_1(\tau^1_Y(y_1,z,\beta))-Y_2(\tau^1_Y(y_2,z,\overline{\beta}^0))\|\le{\cal O}_{T,K}(\|y_1-y_2\|),\\
\displaystyle
\|X_1(\tau^1_X(x_1,w,\overline{\gamma}^0[\beta]))-X_2(\tau^1_X(x_2,w,\gamma[\beta]))\|\le{\cal O}_{T,K}(\|x_1-x_2\|),\\
\displaystyle
\|J_{\overline\tau^1}(x_1,y_1,w,z,\overline\gamma^0[\overline\beta^0],\beta)-J_{\overline\tau^1}(x_2,y_2,w,z,\gamma[\overline\beta^0],\overline\beta^0)\le\\
\displaystyle
\ \ \ {\cal O}_{T,K}(\|(x_1,y_1)-(x_2,y_2)\|),
\end{array}
\end{equation}

\noindent
where $\overline\tau^1=\min\{\tau^1_X(x_2,w,\gamma[\beta]),\tau^1_Y(y_1,z,\beta),T\}$. Our goal is to estimate the difference of the two $J_T$ costs. If $\overline\tau^1=T$, then we are done. Otherwise, we have some cases. We analyze some of them, being the others similarly treated.

1) Suppose that, using the notation of (\ref{eq:unique_sequence}), 

\[
\begin{array}{ll}
\displaystyle
\tau^1=\tau^1_X(x_2,w,\gamma[\beta])<\tau^2=\tau^1_Y(y_1,z,\beta)<\\
\displaystyle
\ \ \ \tau^3=\tau^1_X(x_1,w,\overline{\gamma}^0[\beta])<\tau^4=\tau^1_Y(y_2,z,\overline{\beta}^0),\\
\displaystyle
\tau^2-\tau^1\le{\cal O}_{T,K}(\|(x_1,y_1)-(x_2,y_2)\|).
\end{array}
\]

\noindent
By (\ref{eq:Assumption2-type}), this implies that all four switchings occur in a lap of time not grater than $2{\cal O}_{T,K}(\|(x_1,y_1)-(x_2,y_2)\|)<\mu$, which also implies that, in the meanwhile, no trajectory can switch two times. We then have the pairs $(X_1(\tau^3),-w)$, $X_2(\tau^3),-w)$, as well as the pairs $(Y_1(\tau^4),-z)$, $(Y_2(\tau^4),-z)$.

At the instant $\tau^3$, $X_1$ updates its non-anticipating strategy using $\overline{\gamma}^1$, as given in the view of (\ref{eq:Assumption2}), referring to the points $(X_1(\tau^3),-w)$, $(X_2(\tau^3),-w)$, to the trajectory $X_2(\cdot; X_2(\tau^3),-w,\gamma[\beta](\cdot+\tau^3)])$, and with respect the exit from $\overline\Omega^{-w}_X$, in the time interval $[0,T-\tau^3]$. Similarly, at $\tau^4$, $Y_2$ updates its non-anticipating tuning using $\overline{\beta}^1$ as given in the view of (\ref{eq:Assumption2}), referring to the points $(Y_1(\tau^4),-z)$, $(Y_2(\tau^4),-z)$, to the trajectory $Y_1(\cdot;Y_1(\tau^4),-z,\beta(\cdot+\tau^4))$, and with respect to the exit from $\overline\Omega_Y^{-z}$, in the time interval $[0,T-\tau^4]$. Gluing together, we get the non-anticipating tuning and non-anticipating strategy

\begin{equation}
\label{hatbeta}
\beta\mapsto\hat\beta:t\mapsto
\left\{
\begin{array}{ll}
\displaystyle
\overline{\beta}^0(t), &\mbox{if } 0\le t\le \tau^4,\\
\displaystyle
\overline{\beta}^1(t-\tau^4),&\mbox{if } t>\tau^4,
\end{array}
\right.
\ \ \forall\beta\in{\cal B},
\end{equation}

\begin{equation}
\label{hatgamma}
\beta\mapsto\hat\gamma[\beta]:t\mapsto
\left\{
\begin{array}{ll}
\displaystyle
\overline{\gamma}^0[\beta](t),&\mbox{if } 0\le t\le \tau^3,\\
\displaystyle
\overline{\gamma}^1[\beta](t-\tau^3),&\mbox{if } t>\tau^3,
\end{array}
\right.
\ \ \forall\beta\in{\cal B}.
\end{equation}

Note that \eqref{hatbeta} and \eqref{hatgamma} are non-anticipating in the sense of (\ref{eq:strategies}), (\ref{eq:tuning}), because they use already given non-anticipating constructions and glue them in dependence of the behavior of the trajectories (solutions of (\ref{eq:thermostat_systems})), which are non-anticipating (the state-position only depends on the past behavior).

Let $\tau^5\ge\tau^4$ be a possible subsequent switching instant of one of the trajectories when continuing to move with $\hat\beta$ and $\hat\gamma[\hat\beta]$. Hence, looking to (\ref{eq:sum_costs}), we have

\begin{equation}
\label{eq:splitting}
\begin{array}{ll}
\displaystyle
J_{\tau^5}(x_1,y_1,w,z,\hat\gamma[\hat\beta],\beta)-J_{\tau^5}(x_2,y_2,w,z,\gamma[\hat\beta],\hat\beta)=\\
\displaystyle
J_{\tau^1}(x_1,y_1,w,z,\overline\gamma^0[\overline\beta^0],\beta)-J_{\tau^1}(x_2,y_2,w,z,\gamma[\overline\beta^0],\overline\beta^0)+\\
\displaystyle
\sum_{i=2}^4e^{-\lambda\tau^{i-1}}\left(J^1_{(\tau^i-\tau^{i-1})}-J^2_{(\tau^i-\tau^{i-1})}\right)+\\
\displaystyle
e^{-\lambda{\tau^4}}\Big(J_{(\tau^5-\tau^4)}(X_1(\tau^4),Y_1(\tau^4),-w,-z,\overline\gamma^1[\hat\beta](\cdot+\tau^4-\tau^3),\beta(\cdot+\tau^4))-\\
\displaystyle
J_{(\tau^5-\tau^4)}(X_2(\tau^4),Y_2(\tau^4),-w,-z,\gamma[\hat\beta](\cdot+\tau^4),\overline\beta^1)\Big)
\end{array}
\end{equation}

\noindent
Note that, by semigroup properties, the definition of $\overline\gamma^1$, for $t\ge0$, $\overline\gamma^1[\hat\beta](\cdot+\tau^4-\tau^3)$ corresponds to the non-anticipating strategy constructed for $X_1$ as in (\ref{eq:Assumption2}), taking the points $(X_1(\tau^4),-w)$, $(X_2(\tau^4),-w)$ as references, together with the trajectory $X_2(\cdot;X_2(\tau^4),-w,\gamma[\hat\beta](\cdot+\tau^4))$, and with respect to the exit from $\overline\Omega_X^{-w}$ and the time interval $[\tau^4,T]$. 

By (\ref{eq:Assumption2-type}) and a similar estimate for $J_{(\tau^5-\tau^4)}$, the difference in the second line and the one in the fourth and fifth lines of (\ref{eq:splitting}), are bounded by ${\cal O}_{K,Y}(e^{LT}\|(x_1,y_1)-(x_2,y_2)\|)$, whereas the addenda inside the summation in the third line
are all bounded by $2M{\cal O}_{K,Y}(\|(x_1,y_1)-(x_2,y_2)\|)$ because they consist of integrals of $\ell$ in time interval with length less than ${\cal O}_{K,Y}(\|(x_1,y_1)-(x_2,y_2)\|)$, and $\ell$ is bounded by $M$.

2) Suppose that, using the notation as in the case 1) here above, $\tau^1=\tau^1_X(x_2,w,\gamma[\beta])<\tau^1_Y(y_1,z,\beta)$ but $\tau^1_Y(y_1,z,\beta)-\tau^1_X(x_2,w,\gamma[\beta])>{\cal O}_{T,K}(\|(x_1,y_1)-(x_2,y_2)\|)$. Then it is certainly 

\[
\begin{array}{ll}
\displaystyle
\tau^1=\tau^1_X(x_2,w,\gamma[\beta])\le\tau^2=\tau^1_X(x_1,w,\overline\gamma^0[\beta])<\\
\displaystyle
\tau^3=\tau^1_Y(y_1,z,\beta)\le\tau^4=\tau^1_Y(y_2,z,\overline\beta^0).
\end{array}
\]

\noindent
where $\max\{\tau^2-\tau^1,\tau^4-\tau^3\}\le{\cal O}_{T,K}(\|(x_1,y_1)-(x_2,y_2)\|)$. In this case, $X_1$ update its non-anticipating strategy at $\tau^3$ and $Y_2$ at $\tau^4$, as in the case 1). As in this case, for a subsequent switching instant $\tau^5$, the differences between the costs $J_{(\tau^i-\tau^{i-1})}$ are all estimated in a similar way  as in (\ref{eq:splitting}).

Putting together cases 1),  2) and the others, which are similarly treated, and in particular, considering that, in the time interval $[0,T]$ there can be only a finite number of switching $N_T$, we then  get that the following: 

\begin{Proposition}
\label{prop:almost_continuity}
Given the Main Assumptions and (\ref{eq:controllability-switching}), then, for any $T>0$, for any $K\subset\mathbb{R}^n\times\mathbb{R}^m$ compact there exist $\delta>0$ and a modulus of continuity $\omega_{T,K}$ such that for every $(w,z)$ and for every $(x_1,w_1),(x_2,y_2)\in\overline B_{(w,z)}\cap K$ with $\|(x_1,y_1)-(x_2,y_2)\|\le\delta$, there exist a way to associate $\hat\gamma\in\Gamma$ to any $\gamma\in\Gamma$ and a non-anticipating tuning $\beta\mapsto\hat\beta$ on ${\cal B}$ such that

\begin{equation}
\label{eq:principal-estimate}
J_T(x_1,y_1,w,z,\hat\gamma[\hat\beta],\beta)-J_T(x_2,y_2,w,z,\gamma[\hat\beta],\hat\beta)\le\omega_{T,K}(\|(x_1,y_1)-(x_2,y_2)\|)
\end{equation}

\noindent
Similarly, it holds reversing the roles of $X$ and $Y$, $\alpha\in{\cal A}$ and $\beta\in{\cal B}$, $\gamma\in\Gamma$ and $\xi\in\Xi$.
\end{Proposition}

\section{Continuity and switched DPP}
\label{sec:DPP}
\begin{Proposition}
\label{prop:continuity}
Given the Main Assumptions of Section \ref{sec:assumptions} and \eqref{eq:controllability-switching}, the value functions $\underline{V}$, $\overline V$ \eqref{eq:values} are bounded and continuous in $(x,y)\in\overline B_{(w,z)}$, for all $(w,z)\in\{-1,1\}\times\{-1,1\}$.
\end{Proposition}
{\it Proof:} We prove the proposition for $\underline V$, being the proof for $\overline V$ similar. We are going to use the notations of Section \ref{sec:assumptions}.

The boundedness of $\underline V$ is easily seen, by the boundedness of $\ell$ and the positivity of the discount factor $\lambda>0$.

Let us fix $\varepsilon>0$ and take $T>0$ such that, for all possible trajectories and controls entering the cost $\ell$, it is $\int_T^{+\infty} e^{-\lambda t}\ell dt\le\varepsilon$. Moreover, take a compact $K\subset\mathbb{R}^n\times\mathbb{R}^m$, and, for a fixed pair $(w,z)$ take $(x_1,y_1),(x_2,y_2)\in\overline B_{(w,z)}\cap K$ such that $\|(x_1,y_1)-(x_2,y_2)\|\le\delta$, where $\delta>0$ is given in Proposition \ref{prop:almost_continuity}, with respect to $T$ and $K$.

Take $\gamma_2\in\Gamma$, which realizes $\underline V(x_2,y_2,w,z)$ up to an $\varepsilon$-error, and consider $\hat\gamma_2\in\Gamma$ as the one in Proposition \ref{prop:almost_continuity}, with respect to $T,K$, $(w,z)$, and $(x_1,y_1),(x_2,y_2)$. We get

\begin{equation}
\label{eq:V-V}
\begin{array}{ll}
\displaystyle
\underline V(x_1,y_1,w,z)-\underline V(x_2,y_2,w,z)\le\\
\displaystyle
\sup_{\beta\in{\cal B}}J(x_1,y_1,w,z,\hat\gamma_2[\beta],\beta)-\sup_{\beta\in{\cal B}}J(x_2,y_2,w,z,\gamma_2[\beta],\beta)+\varepsilon\le\\
\displaystyle
J(x_1,y_1,w,z,\hat\gamma_2[\hat\beta_1],\beta_1)-J(x_2,y_2,w,z,\gamma_2[\hat\beta_1],\hat\beta_1)+2\varepsilon,\\
\end{array}
\end{equation}

\noindent
where, $\beta_1\in{\cal B}$ realizes the supremum in the first addendum of the second line up to an $\varepsilon$-error, and $\hat\beta_1$ is as in Proposition \ref{prop:almost_continuity}. Recalling the definition of $T$, using Proposition \ref{prop:almost_continuity}, and continuing with the inequalities (\ref{eq:V-V}), we get 

\[
\begin{array}{ll}
\displaystyle
\underline V(x_1,y_1,w,z)-\underline V(x_2,y_2,w,z)\le\\
\displaystyle
J_T(x_1,y_1,w,z,\hat\gamma_2[\hat\beta_1],\beta_1)-J_T(x_2,y_2,w,z,\gamma_2[\hat\beta_1],\hat\beta_1)+4\varepsilon\le\\
\displaystyle
\omega_{T,K}(\|(x_1,y_1)-(x_2,y_2)\|)+4\varepsilon
\end{array}
\]

\noindent
from which, as usual, by the arbitrariness of $\varepsilon>0$, of the compact $K$ and of the points, we get the required continuity. \cvd

\begin{Proposition}
\label{prop:DPP}
Given the Main Assumptions of Section \ref{sec:assumptions} and \eqref{eq:controllability-switching}, then $\underline V$ and $\overline V$ respectively satisfy

\[
\begin{array}{ll}
\displaystyle
\underline V(x,y,w,z)=\\
\displaystyle
\inf_{\gamma\in\Gamma}\sup_{\beta\in{\cal B}}\left(\int_0^{\tau}e^{-\lambda s}\ell\left(X(s),Y(s),w,z,\gamma[\beta](s),\beta(s)\right)ds+e^{-\lambda\tau}\underline V(X(\tau),Y(\tau),(w,z)^+)\right),\\
\overline V(x,y,w,z)=\\
\displaystyle
\sup_{\xi\in\chi}\sup_{\alpha\in{\cal A}}\left(\int_0^{\tau}e^{-\lambda s}\ell\left(X(s),Y(s),w,z,\alpha(s),\xi[\alpha](s)\right)ds+e^{-\lambda\tau}\overline V(X(\tau),Y(\tau),(w,z)^+)\right),
\end{array}
\]

\noindent
where $X(s)=X(s;x,\gamma[\beta],\beta), Y(s)=Y(s;y,\alpha,\xi[\alpha])$, $\tau$ is the first switching instant and 

\[
(w,z)^+=\left\{
\begin{array}{ll}
\displaystyle
(-w,z)&\mbox{if } \tau=\tau_X<\tau_Y\\
\displaystyle
(-w-z)&\mbox{if } \tau=\tau_X=\tau_Y\\
\displaystyle
(w,-z)&\mbox{if }\ \tau=\tau_Y<\tau_X
\end{array}
\right.
\]

\noindent
is the first ``switched" label.

\end{Proposition}

{\it Proof}. We only prove the equality for $\underline V$. We recall that, in our definition, the switching occurs when the threshold is bypassed. This is the reason for which we consider instants a little bit larger than the switching time, see $\tau_n$ below. We will also use the estimates (\ref{eq:crucial_inequality_thermostat}), which will be discussed in the next section.

Let us denote by $p=(x,y,w,z)\in\overline B_{(w,z)}\times\{(w,z)\}$ any admissible state and by $p_p(\cdot;\gamma[\beta],\beta)$ the corrresponding trajectory. Let us denote by $\omega(p)$ the right-hand side of the equality concerning $\underline{V}$.

Let us fix $\varepsilon>0$ and for any $p'\in\overline B_{(w,z)}\times\{(w,z)\}$, let $\gamma_{p'}\in\Gamma$ be such that

\[
\underline V(p')\ge\sup_{\beta\in{\cal B}}J(p',\gamma_{p'}[\beta],\beta)-\varepsilon.
\]

Claim: $\underline V(p)\le\omega(p)$.

For every $n>0$, $n \in \mathbb{N}$, let us take $\gamma_n\in\Gamma$ such that

\[
\begin{array}{ll}
\displaystyle
\omega_n(p)\ge\sup_{\beta\in{\cal B}}\Big(\int_0^{\tau_n}e^{-\lambda s}\ell\left(X(s),Y(s),w(s),z(s),\gamma_n[\beta](s),\beta(s)\right)ds+\\
\displaystyle
\ \ \ \ e^{-\lambda\tau_n}\underline V(X(\tau_n),Y(\tau_n),w(\tau_n),z(\tau_n)\Big)-\varepsilon,
\end{array}
\]

\noindent
where $\tau_n=\tau_p[\gamma_n,\beta]+1/n$ and

\[
\begin{array}{ll}
\displaystyle
\omega_n(p)=\inf_{\gamma\in\Gamma}\sup_{\beta\in{\cal B}}\Big(\int_0^{\tau_n}e^{-\lambda s}\ell\left(X(s),Y(s),w(s),z(s),\gamma[\beta](s),\beta(s)\right)ds+\\
\displaystyle
\ \ \ \ e^{-\lambda\tau_n}\underline V(X(\tau_n),Y(\tau_n),w(\tau_n),z(\tau_n))\Big).
\end{array}
\]

For $\beta\in{\cal B}$, we define $p_n=p_p(\tau_n;\gamma_n[\beta],\beta)$, and $\delta_n\in\Gamma$ as

\[
\delta_n[\beta](s)=\left\{
\begin{array}{ll}
\displaystyle
\gamma_n[\beta](s),&\mbox{if } 0\le s\le\tau_n,\\
\displaystyle
\gamma_{p_n}[\beta(\cdot+\tau_n)](s-\tau_n), &\mbox{if } s\ge\tau_n.
\end{array}
\right.
\]

Arguing as in Bardi-Capuzzo Dolcetta \cite{barcap} page 437-438, we eventually get

\[
\omega_n(p)\ge\sup_{\beta\in{\cal B}}J(p,\delta_n[\beta],\beta)-2\varepsilon\ge\underline V(p)-2\varepsilon.
\]

By the arbitrariness of $\varepsilon>0$, the claim is proved if we show that $\omega_n(p)-\omega(p)\le{\cal O}(1/n)$ as $n\to+\infty$, where $\cal O$ is an infinitesimal function as its argument tends to zero. For $\varepsilon>0$, take $\gamma_\varepsilon\in\Gamma$ and $\beta_\varepsilon\in{\cal B}$ such that

\begin{equation}
\label{eq:omega_ntoomega}
\begin{array}{ll}
\displaystyle
\omega_n(p)-\omega(p)\le\\
\displaystyle
\int_0^{\tau_n^\varepsilon}e^{-\lambda t}\ell(X(s),Y(s),W(s),Z(s),\gamma_\varepsilon[\beta_\varepsilon](s),\beta_\varepsilon(s))ds+e^{-\lambda\tau_n^\varepsilon}\underline V(p_n)-\\
\displaystyle
\int_0^{\tau^\varepsilon}e^{-\lambda t}\ell(X(s),Y(s),w,z,\gamma_\varepsilon[\beta_\varepsilon](s),\beta_\varepsilon(s)) -e^{-\lambda\tau^\varepsilon}\underline V(p_{switched})+2\varepsilon,
\end{array}
\end{equation}

\noindent
where $\tau^\varepsilon$ is the first switching instant depending on $\gamma_\varepsilon$ and $\beta_\varepsilon$, and $p_{switched}=(X(\tau^\varepsilon),Y(\tau^\varepsilon),(w,z)^+)$. In particular,

\begin{equation}
\label{eq:particular}
\begin{array}{ll}
\displaystyle
\omega(p)\ge\\
\displaystyle
\sup_{\beta\in{\cal B}}\left(\int_0^{\tau^\varepsilon} e^{-\lambda s}\ell\left(X(s),Y(s),w,z,\gamma_\varepsilon[\beta](s),\beta(s)\right)ds
+e^{-\lambda\tau^\varepsilon}\underline V(p_{switched})\right)-\varepsilon.
\end{array}
\end{equation}

\noindent
In order to estimate the second member in (\ref{eq:omega_ntoomega}), we essentially need to compare the values $\underline V(p_n)$ and $\underline V(p_{switched})$, which may have different switching variables, if $p_n$ has an immediate switching after the one of $p_{switch}$ at time $\tau^\varepsilon$.

We denote $p_{switch}=(X(\tau^\varepsilon),Y(\tau^\varepsilon),w,z)$. It is not restrictive to assume that $\gamma_\varepsilon$ is such that, whenever for some $\beta$ it is
$p_{switch}\in\partial\Omega_X^w\times\partial\Omega_Y^z\times\{(w,z)\}$ (i.e. it is a point of possible double switch), $(w,z)^+=(w,-z)$ (that is only $Y$ switches), and
$\underline V(p_{switched})<\underline V(X(\tau^\varepsilon),Y(\tau^\varepsilon),-w,-z)$, then there exists $\zeta>0$ such that the trajectory does not switch (that is $X$ does not switch) in $]\tau^\varepsilon,\tau^\varepsilon+\zeta[$ (if $\zeta$ is sufficiently small, then certainly $Y$ does not switch again, because it has just switched and hence, to do that, it needs to reach the other threshold). Indeed, we can consider $\gamma'_\varepsilon$ defined as

\[
\gamma'_\varepsilon[\beta](t)=\left\{
\begin{array}{ll}
\displaystyle
\gamma_\varepsilon[\beta](t)&\mbox{if } 0\le t\le\tau^\varepsilon\\
\displaystyle
\mbox{and}\\
\displaystyle
\gamma_\varepsilon[\beta](t)&\mbox{if } t\ge\tau^\varepsilon\ \mbox{and } p_{switched}\not\in\partial\Omega_X^w\times\partial\Omega_Y^z\times\{(w,z)\}\\
\displaystyle
\mbox{otherwise:}\\
\displaystyle
\gamma_\varepsilon[\beta](t)&\mbox{if } t\ge\tau^\varepsilon,\ \mbox{and } \underline V(p_{switched})=\underline V(X(\tau^\varepsilon),Y(\tau^\varepsilon),-w,-z)\\
\displaystyle
a_0&\mbox{if } t\ge\tau^\varepsilon,\ \mbox{and } \underline V(p_{switched})<\underline V(X(\tau^\varepsilon),Y(\tau^\varepsilon),-w,-z)\,,
\end{array}
\right.
\]

\noindent
where $a_0$ is inward-pointing in $(X(\tau^\varepsilon),w)$, and we have that $\gamma'_\varepsilon$ still satisfies (\ref{eq:particular}). We can then always assume that if $p_{switch}\in\partial\Omega_X^w\times\partial\Omega^z_Y\times\{(w,z)\}$ then

\begin{equation}
\label{eq:switch}
\begin{array}{ll}
\displaystyle
(w,z)^+=(w,-z)\ \mbox{and } \underline V(p_{switched})<\underline V(X(\tau^\varepsilon),Y(\tau^\varepsilon),-w,-z)\Longrightarrow\\
\displaystyle
\ \ \ (w_n,z_n)=(w,-z)\ \mbox{for large } n,
\end{array}
\end{equation}

\noindent
where $(w_n,z_n)$ is the actual switching variable of $p_n$.
Since $(w_n,z_n)=(w,z)^+$ for large $n$ whenever $p_{switch}\not\in\partial\Omega_X^w\times\partial\Omega^z_Y\times\{(w,z)\}$, as well as whenever $p_{switch}\in\partial\Omega_X^w\times\partial\Omega^z_Y\times\{(w,z)\}$ and  $(w,z)=(-w,-z)$,
in these three cases, by the continuity of $\underline V$, we get the convergence $\underline V(p_n)\to\underline V(p_{switched})$. Two other cases remain when $p_{switch}\in\partial\Omega_X^w\times\partial\Omega^z_Y\times\{(w,z)\}$. In both we use (\ref{eq:crucial_inequality_thermostat}) and the continuity of $\underline V$.

\[
\begin{array}{ll}
\displaystyle
1)\ \ (w,z)^+=(-w,z)\Longrightarrow \underline V(p_{switched})\ge\underline V(p_n)+{\cal O}\left(\frac{1}{n}\right);\\
\displaystyle
2)\ \ (w,z)^+=(w,-z)\ \mbox{and } \underline V(p_{switched})=\underline V(X(\tau^\varepsilon),Y(\tau^\varepsilon),-w,-z)\\
\displaystyle
\Longrightarrow\ \underline V(p_{switched})\ge\underline V(p_n)+{\cal O}\left(\frac{1}{n}\right),
\end{array}
\]

\noindent
where, in 2) we used the fact that, by (\ref{eq:crucial_inequality_thermostat}), it cannot be $\underline V(p_{switched})>\underline V(X(\tau^\varepsilon),Y(\tau^\varepsilon),-w,-z)$, and that, if $(w,z)^+=(w,-z)$ it cannot be $(w_n,z_n)=(-w,z)$ because $z$ is already switched at the time $\tau^\varepsilon$.

\noindent
Hence, in any case we get  $\underline V(p_{switched})\ge\underline V(p_n)+{\cal O}\left(\frac{1}{n}\right)$ and, by the obvious convergence of the integrals in (\ref{eq:omega_ntoomega}), we get the desired estimate.

Claim: $\underline V(x,y,w,z)\ge\omega(x,y,w,z)$.

Arguing as in Bardi-Capuzzo Dolcetta \cite{barcap} page 437--438, we can prove that (with the same notations as in the previous step), for any $\varepsilon>0$

\[
\omega_n(p)\le\underline V(p)+3\varepsilon
\]

The conclusion then still holds because we also have $\omega(p)-\omega_n(p)\le{\cal O}(1/n)$. Indeed, reversing the roles of $\omega$ and $\omega_n$ in (\ref{eq:omega_ntoomega}),(\ref{eq:particular}), we have that if, $(w_n,z_n)=(-w,z)$ or $(w_n,z_n)=(w,-z)$ then $(w,z)^+=(w_n,z_n)$; moreover, if $(w_n,z_n)=(-w,-z)$ for all $n$ then also $(w,z)^+=(-w,-z)$.
\cvd

\section{The HJI systems and uniqueness} 
\label{sec:HJI}

By Proposition \ref{prop:DPP}, for every $(w,z)\in\{-1,1\}$, on $\overline B_{(w,z)}$ the lower value function $\underline V$ of the infinite horizon problem can be interpreted as the lower value function of the exit-time/exit-costs differential game with  dynamics $(x,a)\mapsto f(x,w,a)$, $(y,b)\mapsto g(y,z,b)$, running cost $(x,y,a,b)\mapsto\ell(x,y,w,z,a,b)$ and exit costs (using the same notations as in Section \ref{sec:exit_time})

\[
\Psi_X(\cdot,\cdot)=\underline V(\cdot,\cdot,-w,z),\ \ \Psi_{XY}(\cdot,\cdot)=\underline V(\cdot,\cdot,-w,-z),\ \ \Psi_Y(\cdot,\cdot)=\underline V(\cdot,\cdot,w,-z),
\]

\noindent
and similarly for the upper value function $\overline V$. Moreover, under the controllability hypothesis (\ref{eq:controllability-switching}) such costs satisfy the compatibility hypothesis in (\ref{eq:assumptions}), here stated for $V\in\{\underline V,\overline V\}$: for all $(x,y)\in\partial\Omega_X^w\times\partial\Omega_Y^z$

\begin{equation}
\label{eq:crucial_inequality_thermostat}
V(x,y,w,-z)\le V(x,y,-w,-z)\le V(x,y,-w,z).
\end{equation}

\noindent
Indeed, consider the double switched  state $(x,y,-w,-z)$. Since $X$ minimizes, and since from the $Y$-only switched state $(x,y,w,-z)$ it may freely decide to switch or not, it is

\[
V(x,y,w,-z)\le V(x,y,-w-z).
\]

\noindent
Indeed, from the position $(x,y,w,-z)$ $X$ has at its disposal all the admissible trajectories starting from $(x,y)$ and lying in $\overline B_{(w,-z)}$, at least for a while, as well as all the admissible trajectories starting from $(x,y)$ and immediately moving in
$\overline B_{(-w,-z)}$. By minimization, the previous inequality holds. Symmetrically, it holds for the maximizing player $Y$, getting $V(x,y,-w,-z)\le V(x,y,-w,z)$.

For every fixed $(w,z)$, we define, respectively, the upper Hamiltonian and the lower Hamiltonian for $(x,y,p,q)\in\mathbb{R}^n\times\mathbb{R}^m\times\mathbb{R}^n\times\mathbb{R}^m$, as

\[
\begin{array}{ll}
\displaystyle
UH_{(w,z)}(x,y,p,q)=\min_{b\in B}\max_{a\in A}\{-f(x,w,a)\cdot p-g(y,z,b)\cdot q-\ell(x,y,w,z,a,b)\},\\
\displaystyle
LH_{(w,z)}(x,y,p,q)=\max_{a\in A}\min_{b\in B}\{-f(x,w,a)\cdot p-g(y,z,b)\cdot q-\ell(x,y,w,z,a,b)\}.
\end{array}
\]

By the interpretation as exit-time/exit-costs on every sector $\overline B_{(w,z)}$, by the continuity of the value functions (Proposition \ref{prop:continuity}), by the compatibility (\ref{eq:crucial_inequality_thermostat}), and by Theorem \ref{thm:bagmagzop}, for every fixed $(w,z)$, $\underline V(\cdot,\cdot,w,z)$ and $\overline V(\cdot,\cdot,w,z)$
are the unique bounded and continuous viscosity solutions $u:\overline B_{(w,z)}\to\mathbb{R}$ of the following Isaacs Dirichlet problems in $\overline B_{(w,z)}$ ($=\overline\Omega^w_X\times\overline\Omega^z_Y$), with $H_{(w,z)}=UH_{(w,z)}$ and $H_{(w,z)}=LH_{(w,z)}$, respectively:


\begin{equation}
\label{eq:system_wz}
\left\{
\begin{array}{ll}
\displaystyle
\lambda u(x,y)+H_{(w,z)}(x,y,\nabla_x u(x,y),\nabla_y u(x,y))=0, &\mbox{in } int\overline B_{(w,z)},\\
\displaystyle
u(x,y)=V(x,y,-w,z),&\mbox{on } \partial\Omega_X^w\times\Omega_Y^z,\\
\displaystyle
u(x,y)=V(x,y,-w,-z),&\mbox{on } \partial\Omega_X^w\times\partial\Omega_Y^z,\\
\displaystyle
u(x,y)=V(x,y,w,-z), &\mbox{on } \Omega_X^w\times\partial\Omega_Y^z,
\end{array}
\right.
\end{equation}

\noindent
where the boundary conditions must be also interpreted in the viscosity sense as in (\ref{eq:subsol})--(\ref{eq:supersol}). 

For a function

\[
u:\bigcup_{(w,z)\in\{-1,1\}\times\{-1,1\}}\left(\overline B_{(w,z)}\times(w,z)\right)\to\mathbb{R},
\]

\noindent
we consider the following problem (four Isaacs Dirichlet problems coupled by the boundary conditions in the viscosity sense, that are mutually exchanged)

\begin{equation}
\label{eq:systems}
\left\{
\begin{array}{ll}
\displaystyle
\forall (w,z)\in\{-1,1\}\times\{-1,1\},\ u\ \mbox{solves}\\
\displaystyle
\left\{
\begin{array}{ll}
\displaystyle
\lambda u(x,y,w,z)+H_{(w,z)}(x,y,\nabla_x u,\nabla_y u)=0, &\mbox{in } int\overline B_{(w,z)},\\
\displaystyle
u(x,y,x,w)=u(x,y,-w,z), &\mbox{on } \partial\Omega_X^w\times\Omega_Y^z,\\
\displaystyle
u(x,y,w,z)=u(x,y,-w,-z), &\mbox{on } \partial\Omega_X^w\times\partial\Omega_Y^z,\\
\displaystyle
u(x,y,w,z)=u(x,y,w,-z), &\mbox{on } \Omega_X^w\times\partial\Omega_Y^z.
\end{array}
\right.
\end{array}
\right.
\end{equation}

\begin{Theorem}
\label{thm:fixed-point}
Under the Main Assumptions in Section \ref{sec:assumptions} and the controllability \eqref{eq:controllability-switching}, we get that $\underline V$ (respectively $\overline V$) is the unique bounded and continuous viscosity solution of \eqref{eq:systems} with $H_{(w,z)}=UH_{(w,z)}$ (respectively $H_{(w,z)}=LH_{(w,z)}$), satisfying, for every $(x,y,w,z)\in\partial\Omega_X^w\times\partial\Omega_Y^z\times\{(w,z)\}$

\begin{equation}
\label{eq:boundary_inequality}
u(x,y,w,-z)\le u(x,y,-w,-z)\le u(x,y,-w,z).
\end{equation}

\end{Theorem}

{\it Proof.} 
The fact that $\underline V$ and $\overline V$ are solutions is explained here above. For the uniqueness, 
we are going to use a fixed point argument applied to a suitable functional defined on the subset $C$ of the space

\begin{equation*}
\chi=\left\{u:\bigcup_{w,z\in\{-1,1\}}\left(\overline{B}_{(w,z)}\times\{(w,z)\}\right)\to\mathbb{R}\Big|u\ \mbox{is continuous and bounded}\right\}.
\end{equation*}

\noindent
given by

\[
C=\left\{u\in X\ \Big|\ u\ \mbox{satisfies (\ref{eq:boundary_inequality})} \right\}.
\]

\noindent
Endowed with the uniform convergence topology, $C$ is a complete metric space. Also note that $\underline V,\overline V\in C$. We prove the uniqueness result for the system (\ref{eq:systems}) corresponding to the case $H_{(w,z)}=UH_{(w,z)}$ (i.e. the case solved by $\underline V$). The other case is similar.

The construction of the functional is performed in three steps. 

{\it First step}. We construct a functional $T_1:C\to \chi$ in the following way. Given $\underline u\in C$, for every fixed $(w,z)$, the functions $\underline u(\cdot,\cdot,-w,z)$, $\underline u(\cdot,\cdot,-w,-z)$ and $\underline u(\cdot,\cdot,w,-z)$  give suitable boundary conditions for the sub-problem in (\ref{eq:systems}) with that $(w,z)$ fixed (also compare with (\ref{eq:system_wz})). In particular, they are continuous and satisfy (\ref{eq:boundary_inequality}). Let us denote by $T_1[\underline u;w,z]:\overline B_{(w,z)}\to\mathbb{R}$ such a unique solution. By Theorem \ref{thm:bagmagzop}, $T_1[\underline u;w,z]$ is the lower value function of the exit-time/exit costs differential game with dynamics $(x,a)\mapsto f(x,w,a)$, $(y,b)\mapsto g(y,z,b)$, running cost $(x,y,a,b)\mapsto\ell(x,y,w,z,a,b)$ and exit costs given by the values of $\underline u$ as before. Hence, we define the image of $\underline u\in C$ via $T_1$ as

\begin{equation}
\label{eq:T}
T_1[\underline u]:\bigcup_{(w,z)\in\{-1,1\}\times\{-1,1\}}\left(\overline B_{(w,z)}\times(w,z)\right)\to\mathbb{R},\ \ (x,y,w,z)\mapsto T_1[\underline u;w,z](x,y).
\end{equation}

\noindent
In general, we cannot guarantee that $T_1[\underline u]\in C$, because it may not satisfy (\ref{eq:boundary_inequality}), but certainly it belongs to $\chi$. However, it satisfies similar inequalities as (\ref{eq:boundary_inequality}), that is, for every $(x,y,w,z)\in\partial\Omega_X^w\times\partial\Omega_Y^z\times\{(w,z)\}$

\begin{equation}
\label{eq:boundary_inequality_T}
T_1[\underline u](x,y,w,-z)\le\underline u(x,y,-w,-z)\le T_1[\underline u](x,y,-w,z),
\end{equation}

\noindent
which can be proved similarly to (\ref{eq:crucial_inequality_thermostat}), because, for example, from the point $(x,y,w,-z)$ if $X$ exits (the minimizing player), then the cost paid is $\underline u(x,y,-w,-z)$.

{\it Second step}. We construct a functional $T_2:C\to \chi$ similarly to $T_1$ with the only difference that, the exit costs on the corner points $(x,y,w,z)\in\partial\Omega_X^w\times\partial\Omega_Y^z\times\{(w,z)\}$ are given by $T_1[\underline u](x,y-w,z)$ ($\Psi_X$) and $T_1[\underline u](x,y,w,-z)$ ($\Psi_Y$) if only $X$ or only $Y$ exits, respectively, and by $\underline u(x,y,-w,-z)$ itself ($\Psi_{XY}$) for the case of simultaneous exit. In this way, by (\ref{eq:boundary_inequality_T}), the costs $\Psi_X,\Psi_Y,\Psi_{XY}$ satisfies (\ref{eq:boundary_inequality}). We then construct $T_2[\underline u]$ as in the first step. Again, we have for every $(x,y,w,z)\in\partial\Omega_X^w\times\partial\Omega_Y^z\times\{(w,z)\}$

\[
T_2[\underline u](x,y,w,-z)\le T_1[\underline u](x,y,-w,-z)\le T_2[\underline u](x,y,-w,z).
\]

\noindent
Note that, when $(x,y,w,z)\in\partial\Omega_X^w\times\partial\Omega_Y^z\times\{(w,z)\}$, then from $(x,y,w,-z)$ only $X$ can exit from $\overline B_{(w,-z)}$, and hence, in that case, the paid cost is $T_1[\underline u]$.

{\it Third step}. We construct a functional $T_3:C\to \chi$, as in the second step, but using $T_2[\underline u]$ as exit costs for the exit of $X$ and $Y$ only, and $T_1[\underline u]$ for the simultaneous exit.

It is evident that any solution $u\in C$ of (\ref{eq:systems}) is a fixed point of $T_3$ (as well as of $T_1$ and $T_2$), and in particular $\underline V$ is a fixed point of $T_3$. 
We now prove that $T_3$ is a contraction, and so it admits at most one fixed point, which means that (\ref{eq:systems}) admits at most one solution in $C$, i.e., $\underline V$. From which the proof will be concluded.

Let us take two functions $u^1,u^2\in C$ and a point $(x,y,w,z)$. By construction,
$$
\begin{aligned}
&T_1[u^i](x,y,w,z)=\underline V^i_{(w,z)}(x,y):=\inf_{\gamma\in\Gamma}\sup_{\beta\in{\cal B}}J^i_{(w,z)}(x,y,\gamma[\beta],\beta)=\\
&\inf_{\gamma\in\Gamma}\sup_{\beta\in\mathcal{B}}\Big(\int_0^\tau \ell(X(s),Y(s),\gamma[\beta](s),\beta(s),w,z)+e^{-\lambda\tau}u^i(X(\tau),Y(\tau),(w,z)^+)\Big),
\end{aligned}
$$

\noindent
where $\tau$ is the exit time from $\overline B_{(w,z)}$, and $(w,z)^+$ is $(-w,z),(w,-z), (-w,-z)$ if only $X$ exits, only $Y$ exits, or they simultaneously exit, respectively.

Let us fix $\varepsilon>0$. For suitable 
$\bar{\gamma}\in \Gamma$, $\bar\beta\in{\cal B}$, by definition of infimum and supremum, we have
\begin{equation}\label{T_stima}
	\begin{array}{l}
	\displaystyle
	T_1[u^1](x,y,w,z)-T_1[u^2](x,y,w,z)=\underline V^1_{(w,z)}(x,y)-\underline V^2_{(w,z)}(x,y)\\
	\displaystyle
	\le J^1_{(w,z)}(x,y,\bar\gamma[\bar\beta],\bar\beta)-J^2_{(w,z)}(x,y,\bar\gamma[\bar\beta],\bar\beta)+2\varepsilon.\\
	\end{array}
	\end{equation}	

\noindent
Note that both $J^1$ and $J^2$ are concerning with the same dynamics and running cost, because they are governed by the same controls $\bar\gamma[\bar\beta]$ and $\bar\beta$; starting from the same point $(x,y)$, the possible exit time $\tau$ from $\overline B_{(w,z)}$ is the same, and moreover, the possible switched label $(w,z)^+$ is also the same. Hence, they may differ only for the paid exit cost, which is then paid in the same point $(X(\tau),Y(\tau))$ and with the same discount $e^{-\lambda\tau}$. If the trajectory does not exit, then the difference of the $J^i$'s in the second member of (\ref{T_stima}) is zero. Otherwise, it is of the form $e^{-\lambda\tau}(u^1(X(\tau),Y(\tau),(w,z)^+)-u^2(X(\tau),Y(\tau),(w,z)^+))$. Hence, we have

\begin{equation}
\label{eq:stima_T_1}
\begin{array}{ll}
\displaystyle
T_1[u^1](x,y,w,z)-T_1[u^2](x,y,w,z)\le\\
\displaystyle
	\leq e^{-\lambda \tau}u^1(X(\tau), Y(\tau),(w,z)^+)-u^2(X(\tau), Y(\tau),(w,z)^+)+2\varepsilon.\\
\end{array}
\end{equation}

\noindent
Note that, in the case when the exit occurs at a finite instant of time, if $w^+=-w$ then $X(\tau)$ has a distance from $\partial\Omega_X^{-w}$ equal to $0<\rho_1-\rho_{-1}$, while if $z^+=-z$, $Y(\tau)$ has a distance from $\partial\Omega_Y^{-z}$ equal to $0<\eta_1-\eta_{-1}$.

Applying the same reasoning to $T_2$, we obtain:

\[
\begin{array}{ll}
\displaystyle
T_2[u^1](x,y,w,z)-T_2[u^2](x,y,w,z)\le\\
\displaystyle
e^{-\lambda\sigma}\left(\Psi^1(X(\sigma),Y(\sigma),(w,z)^+)-\Psi^2(X(\sigma),Y(\sigma),(w,z)^+)\right)+2\varepsilon,
\end{array}
\]

\noindent
where $\Psi^i$ is the corresponding exit cost as from the construction of $T_2$ and $\sigma$ is the possible exit time (corresponding to the choice of suitable $\bar\gamma$ and $\bar\beta$ as before). In particular, if only $X$ or only $Y$ exits then $\Psi^i=T_1[u^i]$, otherwise, in case of simultaneous exit, it is $u^i$.

Similarly for $T_3$

\[
\begin{array}{ll}
\displaystyle
T_3[u^1](x,y,w,z)-T_3[u^2](x,y,w,z)\le\\
\displaystyle
e^{-\lambda\nu}\left(\Phi^1(X(\nu),Y(\nu),(w,z)^+)-\Phi^2(X(\nu),Y(\nu),(w,z)^+)\right)+2\varepsilon,
\end{array}
\]

\noindent
where $\Phi^i$ is the corresponding exit cost as from the construction of $T_3$ and $\nu$ is the possible exit time (corresponding to the choice of suitable $\bar\gamma$ and $\bar\beta$ as before). In particular, if only $X$ or only $Y$ exits then $\Phi^i=T_2[u^i]$, otherwise, in case of simultaneous exit, it is $T_1[u^i]$.

Now, suppose that $(X(\nu),Y(\nu),(w,z)^+)=(X(\nu),Y(\nu),-w,z)$, then

\[
\begin{array}{ll}
\displaystyle
T_3[u^1](x,y,w,z)-T_3[u^2](x,y,w,z)\le\\
\displaystyle
e^{-\lambda\nu}\left(T_2[u^1](X(\nu),Y(\nu),-w,z)-T_2[u^2](X(\nu),Y(\nu),-w,z)\right)+2\varepsilon.
\end{array}
\]

\noindent
It is

\[
\begin{array}{ll}
\displaystyle
T_2[u^1](X(\nu),Y(\nu),-w,z)-T_2[u^2](X(\nu),Y(\nu),-w,z)\le\\
\displaystyle
e^{-\lambda\sigma}\left(\Psi^1(X(\sigma),Y(\sigma),(-w,z)^+)-\Psi^2(X(\sigma),Y(\sigma),(-w,z)^+)\right)+2\varepsilon,
\end{array}
\]

\noindent
and suppose that $(-w,z)^+=(w,-z)$, then $\Psi^i=u^i$ and $\|X(\nu)-X(\sigma)\|\ge\rho_1-\rho_{-1}$, which implies

\[
\sigma>\frac{\rho_1-\rho_{-1}}{M}>0,
\]

\noindent
where $M$ is a bound for the dynamics. We then get

\[
\begin{array}{ll}
\displaystyle
T_3[u^1](x,y,w,z)-T_3[u^2](x,y,w,z)\le\\
\displaystyle
e^{-\lambda(\nu+\sigma)}\left(u^1(X(\sigma),Y(\sigma),-w,z)-u^2(X(\sigma),Y(\sigma),-w,z\right)+4\varepsilon\le\\
\displaystyle
e^{-\lambda(\sigma+\nu)}\|u^1-u^2\|_\infty+4\varepsilon.
\end{array}
\]

\noindent
If instead, for example, the sequence of the switching variables is $(w,z)\to(-w,z)\to(-w,-z)\to(-w,z)$, then

\[
\begin{array}{ll}
\displaystyle
T_3[u^1](x,y,w,z)-T_3[u^2](x,y,w,z)\le\\
\displaystyle
e^{-\lambda(\nu+\sigma+\tau)}\left(u^1(X(\tau),Y(\tau),-w,z)-u^2(X(\tau),Y(\tau),-w,z\right)+6\varepsilon\le\\
\displaystyle
e^{-\lambda(\nu+\sigma+\tau)}\|u^1-u^2\|_\infty+6\varepsilon,
\end{array}
\]

\noindent
where 

\[
\tau>\frac{\eta_1-\eta_{-1}}{M}>0,
\]

\noindent
because $\|Y(\sigma)-Y(\tau)\|\ge\eta_1-\eta_{-1}$. Other cases can be
proved in a similar way.

By the arbitrariness of $\varepsilon>0$ and of the point $(x,y,w,z)$, setting 

\[
h=\min\left\{\frac{\rho_1-\rho_{-1}}{M},\frac{\eta_1-\eta_{-1}}{M}\right\}>0,
\]

\noindent
we get

\[
\|T_3[u^1]-T_3[u^2]\|_\infty\le e^{-\lambda h}\|u^1-u^2\|_\infty.
\]

\cvd

\begin{Remark}
\label{rmrk:equilibrium}
By the uniqueness proved in Theorem \ref{thm:fixed-point}, whenever the Hamiltonians $LH$ and $UH$ are equal, then we get $\underline V=\overline V$, that is the differential game has an equilibrium. As usual, the equality of the Hamiltonians can be assured by some particular structure of the running cost $\ell$, for example if, besides the already assumed decoupling feature as in the Main Assumptions in Section \ref{sec:assumptions}, it is also of the form $\ell(x,y,w,z,a,b)=\ell_{(1)}(x,y,w,z)+\ell_{(2)}(a)+\ell_{(3)}(b)$.
\end{Remark}


\section{On numerical treatment}

The numerical treatment of the Isaacs equation is in general a difficult problem, especially in connection with the interpretation of the numerical results in the view of the possible real behavior of the players (see Falcone \cite{fal}, page 494). Moreover, in our case the first main issue would be to obtain numerical results for the exit-time/exit-costs differential game problem in Bagagiolo-Maggistro-Zoppello \cite{bagmagzop}, whose theoretical results are at the basis for the construction of the fixed point procedure applied here to the system of Isaacs equations (\ref{eq:systems}). This seems at the moment a hard question and certainly outside of the goal of the present paper. However, following Falcone \cite{fal}, we write here a possible (certainly not tested) discretization scheme for the present switching problem, under some suitable hypotheses. We just sketch it. 

Notations and hypotheses are as in the previous sections, in particular see (\ref{eq:B}) and the paragraph before (\ref{eq:controllability-switching}). We are going to describe  a possible space-time discretization for the evaluation of the lower value function $\underbar V$. Let us consider two $n$-dimensional and $m$-dimensional compact cubes (centered at the origin and with faces parallel to the axes) $Q_X\subset\mathbb{R}^n$ and $Q_Y\subset\mathbb{R}^m$, such that they contain, in their interior, the switching thresholds $x_1=\rho_1,\rho_{-1}$, $y_1=\eta_1,\eta_{-1}$ on their first axes, respectively. We define $Q=Q_X\times Q_Y$. We assume that $Q$ is invariant for the trajectories, that is, for every $(x,y)\in\partial Q$ and every admissible $(w,z)$, the vectors $(f(x,w,a),g(x,z,b))\in\mathbb{R}^n\times\mathbb{R}^m$ are all inward pointing in $Q$ (starting in $Q$, it is not possible to exit from $Q$). Such an assumption is not in contradiction with the controllability assumption (\ref{eq:controllability-switching}) because the switching boundary is transversal to the boundary of the cube. 

Given a triangle mesh on $Q$ with nodal points $x^i\in Q_X$, $y^j\in Q_Y$, $i,j=1,\dots,K$, and denoting by $h>0$ the step of the time-discretization, we define, for every $(w,z)$, 

\begin{equation}
\label{eq:IandJ}
\begin{array}{ll}
\displaystyle
{\cal I}_w=\left\{i\Big|x^i\in\overline\Omega_X^w\right\},\ \ \ {\cal J}_z=\left\{j\Big|y^j\in\overline\Omega_Y^z\right\},\\
\displaystyle
{\cal I}^{switch}_w=\left\{i\in{\cal I}_w\Big|\exists a\in A\ \mbox{such that } x^i+hf(x^i,w,a)\not\in\overline \Omega_X^w\right\},\\
\displaystyle
{\cal J}^{switch}_z=\left\{j\in{\cal J}_z\Big|\exists b\in B\ \mbox{such that } y^j+hg(y^j,z,b)\not\in\overline \Omega_Y^z\right\},\\
\displaystyle
{\cal I}^{in}_w=\left\{i\in{\cal I}_w\in\Big|i\not\in{\cal I}^{switch}_w\right\},\ \ \ {\cal J}^{in}_z=\left\{j\in{\cal J}_z\in\Big|j\not\in{\cal J}^{switch}_z\right\}.
\end{array}
\end{equation}

For every $w$ (respectively, $z$) and for every $i\in{\cal I}_w$ (respectively, $j\in{\cal J}_z$), we define

\[
A^i_w=\left\{a\in A\Big|x^i+hf(x^i,w,a)\in\overline\Omega_X^w\right\},\ \ \left(\mbox{respectively, } B^j_z=\left\{b\in B\Big|y^j+hg(y^j,z,b)\in\overline\Omega_Y^z\right\}\right).
\]

\noindent
If $h>0$ is small, then, due to the controllability assumptions and the hypotheses on $Q$, both $A^i_w$ end $B^j_z$ are never empty, and, if $i\in{\cal I}_w^{in}$ (respectively, $j\in{\cal J}^{in}_z$) then $A^i_w=A$ (respectively, $B^j_z=B$). Moreover, $(x^i+hf(x^i,w,a),y^j+hg(y^j,z,b))\in Q$, for all $i\in{\cal I}_w$, $j\in{\cal J}_z$, $a\in A$, $b\in B$.

With a suitably constructed mesh, for every $i\in{\cal I}_w$, $a\in A^i_w$ (respectively, $j\in{\cal J}_z$, $b\in B^j_z$), 
we can suitably take convex coefficients $\mu_{i\xi}(a)$ (respectively, $\nu_{j\zeta}(b)$) such that

\begin{equation}
\label{eq:convex}
\sum_{\xi\in{\cal I}_w}\mu_{i\xi}(a)x^\xi=x^i+hf(x^i,w,a),\ \ (\mbox{respectively,}\ \sum_{\zeta\in{\cal J}_z}\nu_{j\zeta}(b)y^\zeta=y^j+hg(y^j,z,b)).
\end{equation}

For every $(w,z)$, let $K_w$ and $K_z$ be the cardinality of ${\cal I}_w$ and ${\cal J}_z$, respectively, and we consider the cartesian product set $\Pi=\bigtimes_{(w,z)\in\{-1,1\}\times\{-1,1\}}\left(\mathbb{R}^{K_w}\times\mathbb{R}^{K_z}\times\{(w,z)\}\right)$, whose elements are denoted by $V:=(V^{(w,z)}):=\left((V^{(w,z)}_{ij})_{i\in{\cal I}_w,j\in{\cal J}_z},w,z\right)_{(w,z)}$. For every $V\in \Pi$ and for every $i\in{\cal I}_w$, $j\in{\cal J}_z$, we define the discretized Isaacs equation

\[
I(V,i,j,w,z)=\max_{b\in B_z^j}\min_{a\in A_w^i}\Big\{\left(1-\lambda h\right)\Big(\sum_{\xi\in{\cal I}_w}\mu_{i\xi}(a)V^{(w,z)}_{\xi j}+\sum_{\zeta\in{\cal J}_z}\nu_{j\zeta}(b)V^{(w,z)}_{i\zeta}\Big)+h\ell(x^i,y^j,w,z,a,b)\Big\}\,,
\]

\noindent
where the coefficients $\mu,\nu$ are defined in (\ref{eq:convex}). We then consider the map

\[
S:\Pi\to \Pi,\ \ \ V\mapsto
S(V):=\left((S^{(w,z)}_{ij}(V))_{i\in{\cal I}_w,j\in{\cal J}_z}\right)_{(w,z)}
\]

\noindent
defined by (componentwise, $(w,z)\in\{-1,1\}\times\{-1,1\}$, $i\in{\cal I}_w$, $j\in{\cal J}_z$)

\begin{equation}
\label{eq:discrete}
S_{ij}^{(w,z)}(V)=
\left\{
\begin{array}{ll}
\displaystyle
I(V,i,j,w,z)&\mbox{if } (i,j)\in{\cal I}^{in}_w\times{\cal J}^{in}_z,\\
\displaystyle
\min\{V_{ij}^{(-w,z)},I(V,i,j,w,z)\}&\mbox{if } (i,j)\in{\cal I}^{switch}_w\times{\cal J}^{in}_z,\\
\displaystyle
\max\{V_{ij}^{(w,-z)},I(V,i,j,w,z)\}&\mbox{if } (i,j)\in{\cal I}^{in}_w\times{\cal J}^{switch}_z,\\
\displaystyle
{\cal V}_2&\mbox{if } (i,j)\in{\cal I}^{switch}_w\times{\cal J}^{switch}_z,\\
\end{array}
\right.
\end{equation}

\noindent
where ${\cal V}_1\le{\cal V}_2\le{\cal V}_3$ is the non-decreasing ordering of the set $\{V^{(w,-z)}_{ij},I(V,i,j,w,z),V^{(-w,z)}_{ij}\}$. 

The lines $2$--$4$ in (\ref{eq:discrete}) represent the discrete version of the boundary conditions in (\ref{eq:system_wz}), also compare with (\ref{eq:lower_dirichlet})--(\ref{eq:supersol}). For $h>0$ small, by the delayed thermostatic switching law and by the boundedness of the dynamics, if $i\in{\cal I}^{switch}_w$ (respectively, $j\in{\cal J}^{switch}_z$) then it is also $i\in{\cal I}^{in}_{-w}$ (respectively, $j\in{\cal J}^{in}_{-z}$). This means that (\ref{eq:discrete}) is well defined. Indeed, if for example $S^{(w,z)}_{ij}(V)$ is defined by the second line in (\ref{eq:discrete}), then the switched component $V^{(-w,z)}_{ij}$ is defined by $I(V,i,j,-w,z)$, which does not involve any other switched components (no "discrete" Zeno phenomenon).  Similarly, if $S^{(w,z)}_{ij}(V)$ is defined by the fourth line of (\ref{eq:discrete}) and, for example, ${\cal V}_2=V^{(w,-z)}_{ij}$, then the definition of ${\cal V}_2$ (i.e. of $V^{(w,-z)}_{ij}$) does not involve another switching in the same variable $z$ (at most two subsequent close swtichings in the two different variables). Finally, we have to restrict the domain of $S$ to the elements $V$ such that $V^{(w,-z)}_{ij}\le V^{(-w,-z)}_{ij}\le V^{(-w,z)}_{ij}$ for all $(w,z),\ i\in{\cal I}_w^{switch},\ j\in{\cal J}^{switch}_z$. 

The function $S$ corresponds to the discretization of the operator $T_1$. We construct $S_2$ and $S_3$ arguing similarly to what was done for the operators $T_2$ and $T_3$. We expect $S_3$ to be a contraction whose unique fixed point is a suitable discretization of the lower value function $\underline V$, evaluated in the nodal points. Then, an interpolated solution can be constructed.
We did not perform qualitative and quantitative studies and numerical tests about the existence and convergence of the discretized solutions. They could be arguments for future studies. 

\section{Conclusions}
In this paper, we have considered an infinite horizon, zero-sum differential game. It is characterized by the fact that the dynamics of each player depend on the evolution of a discrete variable which obey a delayed thermostatic switching law. First, we have proved the continuity of the upper and lower value functions. Second, representing the problem as a coupling of several exit-time differential games, we have characterized each value function as the unique viscosity solution of a system of several Hamilton-Jacobi-Isaacs equations, coupled by the boundary conditions.

The principal hypotheses on the model are a decoupled feature of the dynamics of the two players togehter with a decoupling of the controls in the running cost $\ell$ (see (\ref{eq:thermostat_systems}) and the Main Assumptions below it). The controllability hypothesis (\ref{eq:controllability-switching}) also plays an important role. More general situations are certainly worth studying, in particular, for what concerns the decoupled dynamics and the cost, and they may be the subject of future studies.

The problem here studied is a natural development of the results on exit-time differential games and constrained non-anticipating strategies presented in \cite{bagmagzop} and, up to the knowledge of the authors, is new and original. It could have several applications, for example, in ecological economics for the shallow lake problem or in pursuit evasion game, in which either pursuer or evader dynamics can be affected by a switching discontinuity.

A quantitative application of the present results is surely of interest. As a first step in this direction, in Section 7 we have given some hints and ideas for a possible numerical scheme for our method that could be studied more deeply in future works. 
%

\end{document}